\documentclass{amsart}

\usepackage{xypic, amsmath, amssymb, mathrsfs}

\def\fibdown{\ar@{->>}[d]}
\def\hookdown{\ar@<-.5ex>[d]|{\phantom{a}}|<<{\put(-.7,2){$\scriptstyle\cap$}}}

\newtheorem{thm}{Theorem}[section]
\newtheorem{cor}[thm]{Corollary}
\newtheorem{prop}[thm]{Proposition}
\newtheorem{lem}[thm]{Lemma}

\newenvironment{pf}{\noindent {\it Proof.}}{\hfill$\Box$\\}

\theoremstyle{definition}

\newtheorem{defn}[thm]{Definition}

\newtheorem{re}[thm]{Remark}

\newcommand{\C}{\mathscr{C}}
\newcommand{\D}{\mathscr{D}}

\newcommand{\M}{\mathscr{M}}

\renewcommand{\AA}{\mathbb{A}}

\newcommand{\CC}{\mathbb{C}}

\newcommand{\GG}{\mathbb{G}}

\newcommand{\LL}{\mathbb{L}}

\newcommand{\NN}{\mathbb{N}}

\newcommand{\PP}{\mathbb{P}}

\renewcommand{\SS}{\mathbb{S}}

\newcommand{\ZZ}{\mathbb{Z}}

\newcommand{\lra}{\longrightarrow}

\DeclareMathOperator{\Add}{Add}

\renewcommand{\mod}{{\rm Mod}}

\renewcommand{\i}{\infty}

\DeclareMathOperator{\flag}{Flag}
\renewcommand{\l}{\stackrel}

\DeclareMathOperator{\rings}{Rings}

\DeclareMathOperator{\grass}{Grass}

\makeatletter
\let\c@equation\c@thm
\makeatother
\numberwithin{equation}{section}

\makeatletter
\ifx\SK@label\undefined\let\SK@label\label\fi
 \let\your@thm\@thm
 \def\@thm#1#2#3{\gdef\currthmtype{#3}\your@thm{#1}{#2}{#3}}
 \def\mylabel#1{{\let\your@currentlabel\@currentlabel\def\@currentlabel
  {\currthmtype~\your@currentlabel}
 \SK@label{#1@}}\label{#1}}
 
\makeatother

{\begin{list}{$\bullet$}{\setlength{\labelwidth}{\parindent}\setlength{\leftmargin}{\parindent}}}%
{\end{list}}

\DeclareMathOperator{\spec}{spec}

\DeclareMathOperator{\map}{map}
\DeclareMathOperator{\Map}{Map}

\DeclareMathOperator{\colim}{colim}

\DeclareMathOperator{\im}{im}

\DeclareMathOperator{\sym}{Sym}

\DeclareMathOperator{\holim}{holim}
\DeclareMathOperator{\hocolim}{hocolim}

\renewcommand{\phi}{\varphi}

\newcommand{\eq}[5]{\xymatrix{{#1}\ar[r]&{#2}\ar@<3pt>[r]^-{#3}\ar@<-3pt>[r]_-{#4}&{#5}}}

\DeclareMathOperator{\Grass}{Grass}

\begin{document}

\title[On the Motivic Spectra Representing Cobordism and $K$-Theory]{On the Motivic Spectra Representing\\
Algebraic Cobordism and Algebraic $K$-Theory}

\author{David Gepner and Victor Snaith}


\begin{abstract}
We show that the motivic spectrum representing algebraic $K$-theory is a localization of the suspension spectrum of $\mathbb{P}^\infty$, and similarly that the motivic spectrum representing periodic algebraic cobordism is a localization of the suspension spectrum of $BGL$.
In particular, working over $\mathbb{C}$ and passing to spaces of $\mathbb{C}$-valued points, we obtain new proofs of the topological versions of these theorems, originally due to the second author.
We conclude with a couple of applications:  first, we give a short proof of the motivic Conner-Floyd theorem, and second, we show that algebraic $K$-theory and periodic algebraic cobordism are $E_\infty$ motivic spectra.
\end{abstract}

\maketitle

\section{Introduction}

\subsection{Background and motivation}
Let $(X,\mu)$ be a commutative monoid in the homotopy category of pointed spaces and let $\beta\in\pi_{n}(\Sigma^{\infty} X)$
be an element in the stable homotopy of $X$.
Then $\Sigma^\i X$ is a homotopy commutative ring spectrum, and we may invert the ``multiplication by $\beta$'' map
$$
\mu(\beta):\Sigma^\i X\simeq\Sigma^\i S^0\land\Sigma^\i X\l{\Sigma^{-n}\beta\land 1}{\lra}\Sigma^{-n}\Sigma^\i X\land\Sigma^\i X\l{\Sigma^{-n}\Sigma^\i\mu}{\lra}\Sigma^{-n}\Sigma^\i X.
$$
to obtain a ring spectrum
$$
\Sigma^{\infty} X [1/\beta]:=\colim\{\Sigma^\i X\l{\beta_*}{\lra}\Sigma^{-n}\Sigma^\i X\l{\Sigma^{-n}\beta_*}{\lra}\Sigma^{-2n}\Sigma^\i X\lra\cdots\}
$$
with the property that $\beta_*:\Sigma^\i X[1/\beta]\to\Sigma^{-n}\Sigma^\i X[1/\beta]$ is an equivalence.
In fact, as is well-known, $\Sigma^\i X[1/\beta]$ is universal among $\Sigma^\i X$-algebras $A$ in which $\beta$ becomes a unit.

It was originally shown in \cite{Sn79} (see also \cite{Sn81} for a simpler proof) that the ring spectra $\Sigma^\i_+ BU[1/\beta]$ and $\Sigma^\i_+\CC\PP^\i[1/\beta]$, obtained as above by taking $X$ to be $BU_+$ or $\PP^\i_+$ and $\beta$ a generator of $\pi_2 X$ (a copy of $\ZZ$ in both cases), represent periodic complex cobordism and topological $K$-theory, respectively.
%
%
This motivated an attempt in \cite{Sn79} to define algebraic cobordism by replacing $BGL(\CC)$ in this construction with Quillen's algebraic K-theory spaces \cite{Qu73}.
The result was an algebraic cobordism theory, defined in the ordinary stable homotopy category, which was far too large.

By analogy with topological complex cobordism, algebraic cobordism ought to be the universal oriented algebraic cohomology theory.
However, there are at least two algebraic reformulations of the topological theory; as a result, there are at least two distinct notions of algebraic cobordism popular in the literature today.
One, due to Levine and Morel \cite{LeMo01}, \cite{LeMo02}, constructs a universal ``oriented Borel-Moore'' cohomology theory $\Omega$ by generators and relations in a way reminiscent of the construction of the Lazard ring, and indeed the value of $\Omega$ on the point is the Lazard ring.
However, $\Omega$ is not a generalized motivic cohomology theory in the sense of Morel and Voevodsky \cite{MV99}, so it is not represented by a motivic ring spectrum.

The other notion, and the one relevant to this paper, is Voevodsky's spectrum $MGL$ \cite{Voev98}.
It is a bona fide motivic cohomology theory in the sense that it is defined directly on the level of motivic spectra.
Although the coefficient ring of $MGL$ is still not known (at least in all cases), the orientability of $MGL$ implies that it is an algebra over the Lazard ring, as it carries a formal group law.
Provided one defines an orientation as a compatible family of Thom classes for vector bundles, it is immediate that $MGL$ represents the universal oriented motivic cohomology theory;
moreover, as shown in \cite{PPR07b}, and just as in the classical case, the splitting principle implies that it is enough to specify a Thom class for the universal line bundle.

The infinite Grassmannian
\begin{align*}
BGL_n\simeq\grass_{n,\i}:=\colim_k\grass_{n,k}
\end{align*}
represents, in the $\AA^1$-local homotopy category,
the functor which associates to a variety $X$ the set of isomorphism classes of rank $n$ vector bundles on $X$.
In particular, tensor product of line bundles and Whitney sum of stable vector bundles endow $\PP^\i\simeq BGL_1$ and $BGL\simeq\colim_n BGL_n$ with the structure of abelian group objects in the $\AA^1$-homotopy category.
Note that, over $\CC$, the spaces $\PP^\i(\CC)$ and $BGL(\CC)$ underlying the associated complex-analytic varieties are equivalent to the usual classifying spaces $\CC\PP^\i$ and $BU$.

We might therefore hypothesize, by analogy with topology, that there are equivalences of motivic ring spectra
\begin{align*}
\Sigma^\i_+ BGL[1/\beta]\lra PMGL\qquad\text{and}\qquad
\Sigma^\i_+\PP^\i[1/\beta]\lra K
\end{align*}
where $PMGL$ denotes a periodic version of the algebraic cobordism spectrum $MGL$.
The purpose of this paper is to prove this hypothesis.
In fact, it holds over an arbitrary Noetherian base scheme $S$ of finite Krull dimension, provided one interprets $K$ properly:
the Thomason-Trobaugh $K$-theory of schemes \cite{TT} is not homotopy invariant, and so it cannot possibly define a motivic cohomology theory. Rather, the motivic analogue of $K$-theory is Weibel's homotopy $K$-theory \cite{Weibel88}; the two agree for any {\em regular} scheme.

\subsection{Organization of the paper}
We begin with an overview of the theory of oriented motivic ring spectra.
The notion of an orientation is a powerful one, allowing us to compute first the oriented cohomology of flag varieties and Grassmannians.
We use our calculations to identify the primitive elements in the Hopf algebra $R^0(\ZZ\times BGL)$ with $R^0(BGL_1)$, a key point in our analysis of the abelian group $R^0(K)$ of spectrum maps from $K$ to $R$.

The second section is devoted to algebraic cobordism, in particular the proof that algebraic cobordism is represented by the motivic spectrum $\Sigma^\i_+ BGL[1/\beta]$.
We recall the construction of $MGL$ as well as its periodic version $PMGL$ and note the functors they (co)represent as monoids in the homotopy category of motivic spectra.
We show that $PMGL$ is equivalent to $\bigvee_n\Sigma^\i MGL_n[1/\beta]$ and use the isomorphism $R^0(BGL)\cong\prod_n R^0(MGL_n)$ to identify the functors $\rings(\Sigma^\i_+ BGL[1/\beta],-)$ and $\rings(\bigvee_n\Sigma^\i MGL_n[1/\beta],-)$.

The third section provides the proof that algebraic $K$-theory is represented by the motivic spectrum $\Sigma^\i_+\PP^\i[1/\beta]$.
First we construct a map; to see that it's an equivalence, we note that it's enough to show that the induced map $R^0(K)\to R^0(\Sigma^\i_+\PP^\i[1/\beta])$ is an isomorphism for any $PMGL$-algebra $R$.
An element of $R^0(K)$ amounts to a homotopy class of an infinite loop map $\ZZ\times BGL\simeq\Omega^\i K\to\Omega^\i R$; since loop maps $\ZZ\times BGL\to\Omega^\i R$ are necessarily additive, we are reduced to looking at maps $\PP^\i\to\Omega^\i R$.
We use this to show that the spaces $\map(K,R)$ and $\map(\Sigma^\i_+\PP^\i[1/\beta],R)$ both arise as the homotopy inverse limit of the tower associated to the endomorphism of the space $\map(\Sigma^\i_+\PP^\i,R)$ induced by the action of the Bott map $\PP^1\land\PP^\i\to\PP^\i$, and are therefore homotopy equivalent.

We conclude the paper with a couple of corollaries.
The first is a quick proof of the motivic Conner-Floyd theorem, namely that the map
$$
MGL^{*,*}(X)\otimes_{MGL^{*,*}}K^{*,*}\lra K^{*,*}(X),
$$
induced by an $MGL$-algebra structure on $K$, is an isomorphism for any compact motivic spectrum $X$.
This was first obtained by Panin-Pimenov-R\"ondigs \cite{PPRcf} and follows from a motivic version of the Landweber exact functor theorem \cite{NOS}.
We include a proof because, using the aforementioned structure theorems, we obtain a simplification of the (somewhat similar) method in \cite{PPRcf}, but which is considerably more elementary than that of \cite{NOS}.

Second, it follows immediately from our theorems that both $K$ and $PMGL$ are $E_\i$ {\em as motivic spectra}.
An $E_\i$ motivic spectrum is a {\em coherently} commutative object in an appropriate symmetric monoidal model category of structured motivic spectra, such as P. Hu's motivic $\SS$-modules \cite{Hu} or J.F. Jardine's motivic symmetric spectra \cite{Jard00}; in particular, this is a much stronger than the assertion that algebraic $K$-theory defines a presheaf of (ordinary) $E_\i$ spectra on an appropriate site.
This is already known to be the case for algebraic cobordism, where it is clear from the construction of $MGL$, but does not appear to be known either for {\em periodic} algebraic cobordism or algebraic $K$-theory.

This is important because the category of modules over an $E_\i$ motivic spectrum $R$ inherits a symmetric monoidal structure, at least in the higher categorical sense of \cite{DAGIII}.
As a result, there is a version of derived algebraic geometry which uses $E_\i$ motivic spectra as its basic building blocks.
In \cite{Lur}, J. Lurie shows that $\spec\Sigma^\i_+ \PP^\i[1/\beta]$ is the initial derived scheme over which the derived multiplicative group $\GG_R:=\spec R\land\Sigma^\i_+\ZZ$ acquires an ``orientation'', in the sense that the formal group of $\GG_R$ may be identified with the formal spectrum $\PP^\i\otimes\spec R$.
Since $\Sigma^\i_+\CC\PP^\i[1/\beta]$ represents topological $K$-theory, this is really a theorem about the relation between $K$-theory and the derived multiplicative group, and is the starting point for Lurie's program to similarly relate topological modular forms and derived elliptic curves.
Hence the motivic version of the $K$-theory result may be seen as a small step towards an algebraic version of elliptic cohomology.

%

\subsection{Acknowledgements}
We are very grateful to Mike Hopkins and Rick Jardine for helpful lunch-break discussions during a workshop at the Fields Institute in May 2007.
The first author would also like to thank Sarah Whitehouse for illuminating conversations about operations in $K$-theory, and John Greenlees for his interest in this project as well as its equivariant analogues.

After posting the original version of this paper, we received a communication from P.A. {\O}stv{\ae}r and M. Spitzweck regarding their related project \cite{X}.
The main result of \cite{X} is a different proof that $K\simeq\Sigma^\i_+\PP^\i[1/\beta]$; we thank them for bringing their methods to our attention.

Lastly, we'd like to thank an anonymous referee for a number of constructive comments.

\section{Oriented Cohomology Theories}

\subsection{Motivic spaces}
Throughout this paper, we write $S$ for a Noetherian base scheme of finite Krull dimension.

%
%
%
%
%

\begin{defn}
A {\em motivic space} is a simplicial sheaf on the Nisnevich site of smooth schemes over $S$.
\end{defn}

We often write $0$ for the initial motivic space $\emptyset$, the simplicial sheaf with constant value the set with zero elements, and $1$ for the final motivic space $S$, the simplicial sheaf with constant value the set with one element.

We assume that the reader is familiar with the Morel-Voevodsky $\AA^1$-local model structure on the category of motivic spaces used to define the unstable motivic homotopy category \cite{MV99}.
We adhere to this treatment with one exception:  we adopt a different convention for indexing the simplicial and algebraic spheres.
The {\em simplicial circle} is the pair associated to the constant simplicial sheaves
$$
S^{1,0}:=(\Delta^1,\partial\Delta^1);
$$
its smash powers are the {\em simplicial spheres}
$$
S^{n,0}:=(\Delta^n,\partial\Delta^n).
$$
The {\em algebraic circle} is the multiplicative group scheme $\GG:=\GG_m:=\AA^1-\AA^0$, pointed by the identity section $1\to\GG$; its smash powers define the {\em algebraic spheres}
$$
S^{0,n}:=(\GG,1)^{\land n}.
$$
Putting the two together, we obtain a bi-indexed family of spheres
$$
S^{p,q}:=S^{p,0}\land S^{0,q}.
$$
It is straightforward to show that
$$
(\AA^n-\AA^0,1)\simeq S^{n-1,n}
$$
and
$$
(\AA^n,\AA^n-\AA^0)\simeq(\PP^n,\PP^{n-1})\simeq S^{n,n}.
$$
We emphasize that, according to the more usual grading convention, $S^{p,q}$ is written $S^{p+q,q}$; we find it more intuitive to separate the simplicial and algebraic spheres notationally.
Moreover, for this purposes of this paper, the diagonal spheres
$$S^{n,n}\simeq(\AA^n,\AA^n-\AA^0)\simeq(\PP^n,\PP^{n-1})$$
are far and away the most important, so they will be abbreviated
$$
S^n:=S^{n,n}.
$$
This allows us to get by with just a single index most of the time.

We extend this convention to suspension and loop functors.
That is, $\Sigma(-)$ denotes the endofunctor on pointed motivic spaces (or spectra) defined by
$$
\Sigma X:=S^1\land X:=S^{1,1}\land X.
$$
Similarly, its right adjoint $\Omega(-)$ is defined by
$$
\Omega X:=\map_+(S^1,X):=\map_+(S^{1,1},X).
$$
Note that $\Sigma$ is therefore {\em not} the categorical suspension, which is to say that the cofiber of the unique map $X\to 1$ is given by $S^{1,0}\land X$ instead of $S^{1,1}\land X=\Sigma X$.
While this may be confusing at first, we feel that the notational simplification that results makes it worthwhile in the end.

\subsection{Motivic spectra}
To form the stable motivic category, we formally add desuspensions with respect to the diagonal spheres $S^n=S^{n,n}=(\AA^n,\AA^n-\AA^0)$.

\begin{defn}
A {\em motivic prespectrum} is a sequence of pointed motivic spaces
$$
\{X(0),X(1),\ldots\},
$$
equipped with maps $\Sigma^p X(q)\to X(p+q)$, such that the resulting squares
$$
\xymatrix{
\Sigma^p\Sigma^q X(r)\ar[r]\ar[d] & \Sigma^{p+q} X(r)\ar[d]\\
\Sigma^p X(q+r)\ar[r] & X(p+q+r)}
$$
commute.
\end{defn}

\begin{defn}
A motivic prespectrum is a {\em motivic spectrum} if, for all natural numbers $p,q$, the adjoints $X(q)\to\Omega^{p} X(p+q)$ of the prespectrum structure maps $\Sigma^p X(q)\to X(p+q)$ are weak equivalences.
\end{defn}

A pointed motivic space $X=(X,1)$ gives rise to the suspension spectrum $\Sigma^\i X$, the spectrum associated to the prespectrum with
$$
(\Sigma^\i X)(p):=\Sigma^p X
$$
and structure maps
$$
\Sigma^q\Sigma^p X\lra\Sigma^{p+q} X.
$$
If $X$ isn't already pointed, we usually write $\Sigma^\i_+ X$ for $\Sigma^\i X_+$, where $X_+$ is the pointed space $(X_+,1)\simeq(X,0)$.
If $X$ happens to be the terminal object $1$, we write $\SS:=\Sigma^\i_+ 1$ for the resulting suspension spectrum, the motivic sphere.

We will need that the category of motivic spectra is closed symmetric monoidal with respect to the smash product.
However, we do not focus on the details of its construction, save to say that either P. Hu's theory of motivic $\SS$-modules \cite{Hu} or J.F. Jardine's motivic symmetric spectra \cite{Jard00} will do.

In particular, the category of motivic spectra is tensored and cotensored over itself via the smash product and the motivic function spectrum bifunctors.
We may also regard it as being tensored and cotensored over pointed motivic spaces via the suspension spectrum functor.
Given a motivic spectrum $R$ and a pointed motivic space $X$, we write $X\land R$ for the motivic spectrum $\Sigma^\i X\land R$ and $R^X$ for the motivic spectrum of maps from $\Sigma^\i X$ to $R$.
Here $\Sigma^\i X$ is the motivic spectrum associated to the motivic prespectrum whose value in degree $n$ is the pointed motivic space $\Sigma^n X$.
As a functor from pointed motivic spaces to motivic spectra, $\Sigma^\i$ admits a right adjoint $\Omega^\i$ which associates to a motivic spectrum its underlying motivic ``infinite-loop''.

There are also a number of symmetric monoidal categories over which the category of motivic spectra is naturally enriched.
We write $Y^X$ for the motivic function spectrum of maps from the motivic spectrum $X$ to the motivic spectrum $Y$, $\map(X,Y)=(\Omega^\i Y^X)(S)$ for the (ordinary) space of maps from $X$ to $Y$, and $[X,Y]=Y^0(X)=\pi_0\map(X,Y)$ for the abelian group of homotopy classes of maps from $X$ to $Y$.

\subsection{Motivic ring spectra}
In this paper, unless appropriately qualified, a ring spectrum will always mean a commutative monoid in the homotopy category of motivic spectra.
We reiterate that a motivic spectrum is a $\PP^1$-spectrum; that is, it admits desuspensions by algebraic spheres as well as simplicial spheres.



\begin{defn}
A motivic ring spectrum $R$ is {\em periodic} if the graded ring $\pi_* R$ contains a unit $\mu\in\pi_1 R$ in degree one.
\end{defn}

\begin{re}
Since $\pi_1 R$ is by definition $\pi_0\map_+(\PP^1,\Omega^\i R)$, and over $\spec\CC$, $\PP^1(\CC)\simeq\CC\PP^1$, the topological $2$-sphere, this is compatible with the notion of an {\em even periodic} ring spectrum so common in ordinary stable homotopy theory.
\end{re}

\begin{prop}
If $R$ is periodic then $R\simeq\Sigma^n R$ for all $n$.
\end{prop}

\begin{pf}
Let $\mu\in\pi_1 R$ be a unit with inverse $\mu^{-1}\in\pi_{-1}R$.
Then for any $n$, the multiplication by $\mu^{-n}$ map
$$
R\lra\Sigma^{n} R
$$
is an equivalence, since multiplication by $\mu^n$ provides an inverse.
\end{pf}


Let $P\SS$ denote the periodic sphere, the motivic spectrum
$$
P\SS:=\bigvee_{n\in\ZZ}\Sigma^n\SS.
$$
With respect to the multiplication induced by the equivalences $\Sigma^p\SS\land\Sigma^q\SS\to\Sigma^{p+q}\SS$, the unit in degree one given by the inclusion $\Sigma^1\SS\to P\SS$ makes $P\SS$ into a periodic $\SS$-algebra.
More generally, given an arbitrary motivic ring spectrum $R$,
$$
PR:=P\SS{\land}R\simeq\bigvee_{n\in\ZZ}\Sigma^n R
$$
is a periodic ring spectrum equipped with a ring map $R\to PR$.

\begin{prop}
Let $R$ be a motivic ring spectrum.
Then homotopy classes of ring maps $P\SS\to R$ naturally biject with units in $\pi_1 R$.
\end{prop}

\begin{pf}
By definition, ring maps $P\SS\to R$ are indexed by families of elements $r_n\in\pi_n R$ with $r_m r_n=r_{m+n}$ and $r_0=1$.
Hence $r_n=r_1^n$, and in particular $r_{-1}=r_1^{-1}$.
\end{pf}

Said differently, the homotopy category of periodic motivic ring spectra is equivalent to the full subcategory of the homotopy category of motivic ring spectra which admit a ring map from $P\SS$.
This is not the same as the homotopy category of $P\SS$-algebras, in which only those maps which preserve the distinguished unit are allowed.

\begin{cor}
Let $Q$ be a motivic ring spectrum and $R$ a periodic motivic ring spectrum.
Then the set of homotopy classes of ring maps $PQ\to R$ is naturally isomorphic to the set of pairs consisting of a homotopy class of ring map $Q\to R$ and a distinguished unit $\mu\in\pi_1 R$.
\end{cor}


%

\subsection{Orientations}
Let $R$ be a motivic ring spectrum.

\begin{defn}
The {\em Thom space} of an $n$-plane bundle $V\to X$ is the pair $(V,V-X)$, where $V-X$ denotes the complement in $V$ of the zero section $X\to V$.
\end{defn}

Given two vector bundles $V\to X$ and $W\to Y$, the Thom space $(V\times W,V\times W-X\times Y)$ of the product bundle $V\times W\to X\times Y$ is equivalent (even isomorphic) to the smash product $(V,V-X)\land (W,W-Y)$ of the Thom spaces.
Since the Thom space of the trivial $1$-dimensional bundle $\AA^1\to\AA^0$ is the motivic $1$-sphere $S^1\simeq(\AA^1,\AA^1-\AA^0)$, we see that the Thom space of the trivial $n$-dimensional bundle $\AA^n\to\AA^0$ is the motivic $n$-sphere $S^n\simeq(\AA^n,\AA^n-\AA^0)$.
Note that the complement of the zero section $\LL-\PP^\i$ of the universal line bundle $\LL\to\PP^\i\simeq B\GG$ is equivalent to the total space of the universal principal $\GG$-bundle $E\GG\to B\GG$, which is contractible.
Hence the Thom space of $\LL\to\PP^\i$ is equivalent to $(\PP^\i,\PP^0)$, and the Thom space of the restriction of $\LL\to\PP^\i$ along the inclusion $\PP^1\to\PP^\i$ is equivalent to $(\PP^1,\PP^0)\simeq S^1$.

\begin{defn}
An {\em orientation} of $R$ is the assignment, to each $m$-plane bundle $V\to X$, of a class $\theta(V/X)\in R^m(V,V-X)$, in such a way that\smallskip
\begin{enumerate}
\item
for any $f:Y\to X$, the class $\theta(f^*V/Y)$ of the restriction $f^*V\to Y$ of $V\to X$ is equal to the restriction $f^*\theta(V/X)$ of the class $\theta(V/X)$ in $R^m(f^*V,f^*V-Y)$,\smallskip
\item
for any $n$-plane bundle $W\to Y$, the (external) product $\theta(V/X)\times\theta(W/Y)$ of the classes $\theta(V/X)$ and $\theta(W/Y)$ is equal to the class $\theta(V\times W/X\times Y)$ of the (external) product of $V\to X$ and $W\to Y$ in $R^{m+n}(V\times W,V\times W-X\times Y)$, and\smallskip
\item
if $\LL\to\PP^\i$ is the universal line bundle and $i:\PP^1\to\PP^\i$ denotes the inclusion, then $i^*\theta(\LL/\PP^\i)\in R^1(f^*\LL,f^*\LL-\PP^1)$ corresponds to $1\in R^0(S^0)$ via the isomorphism $R^0(S^0)\cong R^1(S^1)\cong R^1(f^*\LL,f^*\LL-\PP^1)$.
\end{enumerate}
\end{defn}

Given an orientation of $R$, the class $\theta(V/X)\in R^n(V,V-X)$ associated to a $n$-plane bundle $V\to X$ is called the {\em Thom class} of $V\to X$.
The main utility of Thom classes is that they define $R^*(X)$-module isomorphisms $R^*(X)\to R^{*+n}(V,V-X)$ (cf. \cite{PPR07b}).

\begin{re}
The naturality condition implies that it is enough to specify Thom classes for the universal vector bundles $V_n\to BGL_n$.
We write $MGL_n$ for the Thom space of $V_n\to BGL_n$ and $\theta_n$ for $\theta(V_n/BGL_n)\in R^n(MGL_n)$.
\end{re}

%

\subsection{Basic calculations in oriented cohomology}
In this section we fix an oriented motivic ring spectrum $R$ equipped with a unit $\mu\in\pi_1 R$.
Note that we can use $\mu$ to move the Thom classes $\theta_n\in R^n(MGL_n)$ to degree zero Thom classes $\vartheta_n:=\mu^n\theta_n\in R^0(MGL_n)$.
The following calculations are well known (cf. \cite{Ad74}, \cite{At68}, \cite{PPR07b}).
Note that all (co)homology is implicitly the (co)homology of a pair.
In particular, if $X$ is unpointed, then $R^0(X):=R^0(X,0)$, where $0\to X$ is the unique map from the inital object $0$; if $X$ is pointed, then $R^0(X):=R^0(X,1)$, where $1\to X$ is the designated map from the terminal object $1$.



\begin{prop}\label{RcohPn}
The first Chern class of the tautological line bundle on $\PP^n$ defines a ring isomorphism $R^0[\lambda]/(\lambda^{n+1})\to R^0(\PP^n)$.
\end{prop}

\begin{pf}
Inductively, one has a morphism of exact sequences
$$
\xymatrix{
\lambda^n R^0[\lambda]/(\lambda^{n+1})\ar[r]\ar[d] & R^0[\lambda]/(\lambda^{n+1})\ar[r]\ar[d] & R^0[\lambda]/(\lambda^n)\ar[d]\\
R^0(\PP^n,\PP^{n-1})\ar[r]                   & R^0(\PP^n)\ar[r]           & R^0(\PP^{n-1})}
$$
in which the left and right, and hence also the middle, vertical maps are isomorphisms.
\end{pf}

\begin{prop}\label{RcohPi}
The first Chern class of the tautological line bundle on $\PP^\i$ defines a ring isomorphism $R^0[\![\lambda]\!]\cong\lim R^0[\lambda]/(\lambda^n)\to R^0(\PP^\i_+)$.
\end{prop}

\begin{pf}
The $\lim^1$ term in the exact sequence
$$
0\lra{\lim}^1 R^{-1,0}(\PP^n)\lra R^{0}(\PP^\i)\lra\lim R^{0}(\PP^n)
$$
vanishes because the maps
$R^{-1,0}(\PP^n)\lra R^{-1,0}(\PP^{n-1})$
are surjective.
\end{pf}

\begin{cor}
For each $n$, the natural map
$$
R_0(\PP^n)\lra\hom_{R^0}(R^0(\PP^n),R^0)
$$
is an isomorphism.
\end{cor}

\begin{pf}
The dual of \ref{RcohPn} shows that $R_0(\PP^n)$ is free of rank $n+1$ over $R_0$.
\end{pf}

%


\begin{prop}[Atiyah \cite{At68}]\label{lh}
Let $p:Y\to X$ be a map of quasicompact $S$-schemes and let $y_1,\ldots,y_n$ be elements of $R^{0}(Y)$.
Let $M$ be the free abelian group on the $y_1,\ldots,y_n$, and suppose that $X$ has a cover by open subschemes $U$ such that for all open $V$ in $U$, the natural map
$$
R^{0}(V)\otimes M\lra R^{0}(p^{-1}V)
$$
is an isomorphism.
Then, for any open $W$ in $X$, the map
$$
R^{0}(X,W)\otimes M\lra R^{0}(Y,p^{-1}W)
$$
is an isomorphism.
\end{prop}

\begin{pf}
Apply Atiyah's proof \cite{At68}, {\em mutatis mutandis}.
\end{pf}

\begin{prop}\label{trivkun}
Let $Z$ be an $S$-scheme such that, for any homotopy commutative $R$-algebra $A$, $A^0(Z)\cong R^0(Z)\otimes_{R^0}A^0$.
Then, for any $S$-scheme $X$, $R^0(Z\times X)\cong R^0(Z)\otimes_{R^0} R^0(X)$.
\end{prop}

\begin{pf}
The diagonal of $X$ induces a homotopy commutative $R$-algebra structure on $A=R^X$, the cotensor of the motivic space $X$ with the motivic spectrum $R$.
Hence
$$A^0(Z)\cong R^0(Z)\otimes_{R^0}A^0\cong R^0(Z)\otimes_{R^0} R^0(X).$$
\end{pf}

\begin{cor}
Let $p:V\to X$ be a rank $n$ vector bundle over a quasicompact $S$-scheme $X$ and let $L\to\PP(V)$ be the tautological line bundle.
Then the map which sends $\lambda$ to the first Chern class of $L$ induces an isomorphism
$$
R^{0}(X)[\lambda]/(\lambda^n-\lambda^{n-1}c_1 V+\cdots +(-1)^n c_n V)\longrightarrow R^{0}(\PP(V))
$$
of $R^{0}$-algebras.
\end{cor}

\begin{pf}
If $V$ is trivial then $\PP(V)\cong\PP^{n-1}_X$, and the result follows from Propositions \ref{RcohPn} and \ref{trivkun}.
In general, the projection $\PP(V)\to X$ is still locally trivial, so we may apply Proposition \ref{lh}, with $W$ empty and $\{y_i\}$ the image in $R^0(\PP(V))$ of a basis for $R^0(X)[\lambda]/(\lambda^n-\ldots+(-1)^n c_n V)$ as a free $R^0(X)$-module.
\end{pf}


\begin{prop}\label{flag}
Let $V\to X$ be a rank $n$ vector bundle over a quasicompact $S$-scheme $X$, $\flag(V)\to X$ the associated flag bundle, and $\sigma_k(x_1,\ldots,x_n)$, $1\leq k\leq n$, the $k^{\text th}$ elementary symmetric function in the indeterminates $\lambda_i$.
Then the map
$$
R^{0}(X)[\lambda_1,\ldots,\lambda_n]/(\{c_k(V)-\sigma_k(\lambda_1,\ldots,\lambda_n)\}_{k>0})\lra R^{0}(\flag(V))
$$
which sends the $\lambda_i$ to the first Chern classes of the $n$ tautological line bundles on $\flag(V)$, is an isomorphism of $R^0$-algebras.
\end{prop}

\begin{pf}
The evident relations among the Chern classes imply that the map is well-defined.
Using Proposition \ref{lh} and a basis for the free $R^0$-module $R^0(\flag(\AA^{n-1})$, it follows inductively from the fibration $\flag(\AA^{n-1})\to\flag(\AA^n)\to\PP^{n-1}$ that
$$
R^0[\lambda_1,\ldots,\lambda_n]/(\{\sigma_k(\lambda_1,\ldots,\lambda_n)\})\lra R^0(\flag(\AA^n))
$$
is an isomorphism.
Using Proposition \ref{trivkun}, we deduce the desired result for trivial vector bundles $V$.
For the general case, we apply Proposition \ref{lh} again, with a basis of the free $R^0(X)$-module $R^0(X)[\lambda_1,\ldots,\lambda_n]/(\{c_k(V)-\sigma_k\})$ giving the necessary elements of $R^0(\flag(V))$.
\end{pf}

\begin{prop}\label{cograss}
Let $p:V\to X$ be an rank $n$ vector bundle over a quasicompact $S$-scheme $X$, let $q:\grass_m(V)\to X$ be the Grassmannian bundle of $m$-dimensional subspaces of $V$, let $\xi_m(V)\to\grass_m(V)$ be the tautological $m$-plane bundle over $\grass_m(V)$, and write $q^*(V)/\xi_m(V)$ for the quotient $(n-m)$-plane bundle.
Then the map
$$
R^{0}(X)[\sigma_1,\ldots,\sigma_m,\tau_1,\ldots,\tau_{n-m}]\longrightarrow R^{0}(\grass_m(V))
$$
which sends $\sigma_i$ to $c_i(\xi_m(V))$ and $\tau_j$ to $c_j(q^*(V)/\xi_m(V))$ induces an isomorphism
$$
R^{0}(X)[\sigma_1,\ldots,\sigma_m,\tau_1,\ldots,\tau_{n-m}]/(\{c_k(V)-\underset{i+j=k}\Sigma\sigma_i\tau_j\})\lra R^{0}(\grass_m(V))
$$
of $R^{0}$-algebras (as usual, $c_k(V)=0$ for $k>n$ and $c_0(V)=\sigma_0=\tau_0=1$).
\end{prop}

\begin{pf}
The identity $q^*c(V)=c(q^*V)=c(\xi_m(V))c(q^*V/\xi_m(V))$ implies that each $c_k(V)-\Sigma\sigma_i\tau_j$ is sent to zero, so the map is well-defined.
Just as in the case of flag bundles, use induction together with the fibration $\flag(\AA^m)\times\flag(\AA^{n-m})\to\flag(\AA^n)\to\grass_m(\AA^n)$ to see that
$$R^0[\sigma_1,\ldots,\sigma_m,\tau_1,\ldots,\tau_{n-m}]/(\{\underset{i+j=k}\Sigma\sigma_i\tau_j\})\longrightarrow R^0(\grass_m(\AA^n))$$
is an isomorphism.
This implies the result for trivial vector bundles by Proposition \ref{trivkun}, and we deduce the general result from Proposition \ref{lh}.
\end{pf}


%

%
%
%
%
%
%

\begin{prop}\label{RcohBGLn}
There are isomorphisms
$$
R^{0}(BGL_n)\lra R^{0}(BGL_1^n)^{\Sigma_n}\cong R^0[\![\lambda_1,\ldots,\lambda_n]\!]^{\Sigma_n}\cong R^0[\![\sigma_1,\ldots,\sigma_n]\!].
$$
\end{prop}

\begin{pf}
Writing $V_n\to BGL_n$ for the tautological vector bundle, we have an equivalence $\flag(V_n)\simeq BGL_1^n$.
Inductively, we have isomorphisms
$$
R^{0}[\![\lambda_1,\ldots,\lambda_n]\!]\lra R^{0}(BGL_1^n)
$$
and the map
$$
R^{0}(BGL_n)\lra R^{0}(BGL_1^n)\cong R^{0}[\![\lambda_1,\ldots,\lambda_n]\!]
$$
factors through the invariant subring $R^{0}[\![\lambda_1,\ldots,\lambda_n]\!]^{\Sigma_n}\cong R^0[\![\sigma_1,\ldots,\sigma_n]\!]$.
By Proposition \ref{flag}, $R^{0}(BGL_1^n)$ is free of rank $n!$ over $R^{0}(BGL_n)$, so it follows that $R^{0}(BGL_n)\cong R^{0}(BGL_1^n)^{\Sigma_n}$.
\end{pf}

\begin{cor}\label{cohdualho}
The natural map
$$
R^0(BGL_n)\lra\hom_{R_0}(\sym^n_{R_0} R_0(\PP^\i), R_0)
$$
is an isomorphism.
\end{cor}

\begin{pf}
By Proposition \ref{RcohBGLn}, we need only check this for $n=1$.
But
$$R^0(\PP^m)\cong\hom_{R_0}(R_0(\PP^m),R_0),$$
both being free of rank $m+1$ over $R^0$, and
$$
R^0(\PP^\i)\cong\lim R^0(\PP^m)\cong\hom_{R_0}(\colim R_0(\PP^m),R_0)\cong\hom_{R_0}(R_0(\PP^\i),R_0)
$$
by Proposition \ref{RcohPn}.
\end{pf}

%

\begin{cor}
There are isomorphisms $R^{0}(BGL)\cong\lim_n R^{0}(BGL_n)\cong R^{0}[\![\sigma_1,\sigma_2,\ldots]\!]$.
\end{cor}

\begin{pf}
The $\lim^1$ term in the short exact sequence
$$
0\lra{\lim}^1_n R^{1,0}(BGL_n)\lra R^{0}(BGL)\lra{\lim}_n R^{0}(BGL_n)
$$
vanishes since the maps $R^{1,0}(BGL_n)\to R^{1,0}(BGL_{n-1})$ are surjective.
\end{pf}

\subsection{The oriented cohomology of $BGL_+\land Z$}
Let $R$ be an oriented periodic motivic ring spectrum and let $Z$ be an arbitrary motivic spectrum.
Recall (cf. \cite{DI}) that a motivic spectrum is {\em cellular} if belongs to the smallest full subcategory of motivic spectra which is closed under homotopy colimits and contains the spheres $S^{p,q}$ for all $p,q\in\ZZ$, and that a motivic space is {\em stably cellular} if its suspension spectrum is cellular.

\begin{prop}\label{Kunneth}
Let $X:=\colim_n X_n$ be a telescope of finite cellular motivic spectra such that each $R^{*,*}(X_n)$ is a finite free $R^{*,*}$-module and, for any motivic spectrum $Z$,
the induced maps
$$
R^0(X_n\land Z)\longrightarrow R^0(X_{n-1}\land Z)
$$
are surjective.
Then the natural map
$$
R^0(X)\underset{R^0}{\widehat{\otimes}}R^0(Z)\lra R^0(X\land Z)
$$
is an isomorphism.
\end{prop}

\begin{pf}
This is an immediate consequence of the motivic K\"unneth spectral sequence of Dugger-Isaksen \cite{DI}.
Indeed, for each $n$, $R^{*,*}(X_n)$ is a free $R^{*,*}$-module, so the spectral sequence
$$
\text{Tor}_*^{R^{*,*}}(R^{*,*}(X_n),R^{*,*}(Z))\Rightarrow R^{*,*}(X_n\land Z)
$$
collapses to yield the isomorphism
$$
R^{*,*}(X_n)\underset{R^{*,*}}{\widehat{\otimes}}R^{*,*}(Z)\cong R^{*,*}(X_n\land Z).
$$
Moreover, by hypothesis, each of the relevant $\lim^1$ terms vanish, so that
\begin{align*}
R^{*,*}(X)\underset{R^{*,*}}{\widehat{\otimes}} R^{*,*}(Z)&\cong\lim_n R^{*,*}(X_n)\underset{R^{*,*}}{{\otimes}} R^{*,*}(Z)\\
&\cong\lim_n R^{*,*}(X_n\land Z)\cong R^{*,*}(X\land Z).
\end{align*}
\vspace{-.5cm}
\end{pf}

\begin{cor}\label{cohomologyPtensorZ}
Let $Z$ be a motivic spectrum.
Then there are natural isomorphisms
$$
R^0(\PP^\i)\underset{R^0}{\widehat{\otimes}}R^0(Z)\lra R^0(\PP^\i_+\land Z)
$$
and
$$
R^0(BGL)\underset{R^0}{\widehat{\otimes}}R^0(Z)\lra R^0(BGL_+\land Z)
$$
\end{cor}

\begin{pf}
For each $m$,
$$
BGL_m\simeq\colim_n\Grass_m(\AA^n)
$$
is a colimit of finite stably cellular motivic spaces such that, for each $n$, $R^{*,*}(\Grass_m(\AA^n))$ is a free $R^{*,*}$-module and
\begin{eqnarray*}
R^0(\Grass_m(\AA^n))\otimes_{R^0} R^0(Z)&\cong &R^0(\Grass_m(\AA^n)_+\land Z)\lra\\
R^0(\Grass_m(\AA^{n-1})_+\land Z)&\cong &R^0(\Grass_m(\AA^{n-1}))\otimes_{R_0} R^0(Z)
\end{eqnarray*}
is (split) surjective.
It therefore follows from Proposition \ref{Kunneth} that, for each $m$,
$$
R^0(\Grass_m(\AA^\i)_+)\underset{R^0}{\widehat{\otimes}} R^0(Z)\cong R^0(\Grass_m(\AA^\i)_+\land Z).
$$
Taking $m=1$ yields the result for $\PP^\i$; for $BGL$, we must consider the sequence
$$
BGL\simeq\colim_m BGL_m\simeq\Grass_m(\AA^\i)
$$
in which the maps come from a fixed isomorphism $\AA^1\oplus\AA^\i\cong\AA^\i$.
Note, however, that it follows from the above, together with the (split) surjection
$$
R^0(\Grass_m(\AA^\i)\lra R^0(\Grass_{m-1}(\AA^\i))
$$
of Proposition \ref{RcohBGLn}, that, for each $m$,
$$
R^0(\Grass_m(\AA^\i)_+\land Z)\lra R^0(\Grass_{m-1}(\AA^\i)_+\land Z)
$$
is (split) surjective, so that the $\lim^1$ term vanishes and
$$
R^0(BGL)\underset{R^0}{\widehat{\otimes}} R^0(Z)\cong\lim_m R^0(\Grass_m(\AA^\i)_+\land Z)\cong R^0(BGL_+\land Z).
$$
\end{pf}

\subsection{Primitives in the oriented cohomology of $BGL$}
Let $R$ be an oriented periodic motivic ring spectrum.
As is shown in Section 4.3 of \cite{MV99}, the group completion
$$BGL_\ZZ\simeq\Omega^{1,0}B(BGL_\NN)$$
(usually written $\ZZ\times BGL$) of the additive monoid $BGL_\NN=\coprod_{n\in\NN} BGL_n$ fits into a fibration sequence
$$
BGL\lra BGL_\ZZ\lra\ZZ,
$$
where $BGL\simeq\colim_n BGL_n$.
As $BGL_\NN$ is commutative up to homotopy, $BGL_\ZZ$ is an abelian group object in the motivic homotopy category.


\begin{lem}\label{add}
Let $\Add(BGL_\ZZ,\Omega^\i R)$ denote the abelian group of homotopy classes of additive maps $BGL_\ZZ\to\Omega^\i R$.
Then the inclusion
$$
\Add(BGL_\ZZ,\Omega^\i R)\lra R^{0}(BGL_\ZZ)
$$
identifies $\Add(BGL_\ZZ,\Omega^\i R)$ with the abelian group of primitive elements in the Hopf algebra $R^{0}(BGL_\ZZ)$.
\end{lem}

\begin{pf}
By definition, there is an equalizer diagram
$$
\xymatrix{
\Add(BGL_\ZZ,\Omega^\i R)\ar[r] & R^{0}(BGL_\ZZ)\ar@<2pt>[r]\ar@<-2pt>[r] & (BGL_\ZZ\times BGL_\ZZ)}
$$
associated to the square
$$
\xymatrix{
BGL_\ZZ\times BGL_\ZZ\ar[r]\ar[d] & BGL_\ZZ\ar[d]\\
\Omega^\i R\times\Omega^\i R\ar[r] & \Omega^\i R}
$$
in which the horizontal maps are the addition maps.
Let $\delta$ denote the Hopf algebra diagonal
$$
\delta:R^0(BGL_\ZZ)\lra R^{0}(BGL_\ZZ\times BGL_\ZZ)\cong R^{0}(BGL_\ZZ)\underset{R^{0}}{\widehat{\otimes}}R^{0}(BGL_\ZZ).
$$
Then the equalizer consists of those $f\in R^0(BGL_\ZZ)$ such that
$\delta(f)=f\otimes 1+1\otimes f$.
This identifies $\Add(BGL_\ZZ,\Omega^\i R)$ with the primitive elements in $R^{0}(BGL_\ZZ)$.
\end{pf}

\begin{lem}\label{redadd}
There are natural isomorphisms
\begin{align*}
\Add(BGL_\ZZ,\Omega^\i R)&\cong\Add(BGL,\Omega^\i R)\times\Add(\ZZ,\Omega^\i R)\\
&\cong\Add(BGL,\Omega^\i R)\times R^0.
\end{align*}
\end{lem}

\begin{pf}
The product of additive maps is additive, and, in any category with finite products and countable coproducts, $\ZZ=\coprod_{\ZZ} 1$ is the free abelian group on the terminal object $1$.
\end{pf}

\begin{prop}\label{split}
The map
$$\Add(BGL_\ZZ,\Omega^\i R)\lra R^{0}(BGL_1),$$
obtained by restricting an additive map $BGL_\ZZ\to\Omega^\i R$ along the inclusion $BGL_1\to BGL_\ZZ$, is an isomorphism.
\end{prop}

\begin{pf}
By Lemma \ref{redadd}, it's enough to show that the inclusion $(BGL_1,1)\to (BGL,1)$ induces an isomorphism
$$
\Add(BGL,\Omega^\i R)\lra R^0(BGL_1,1).
$$
Thus let $M=R_0(BGL_1,1)$, and consider the $R_0$-algebra
$$
A:=\bigoplus_{n\geq 0}\sym^n_{R_0} M
$$
together with its augmentation ideal
$$
I:=\bigoplus_{n>0}\sym^n_{R_0} M.
$$
We have isomorphisms of split short exact sequences
$$
\xymatrix{
0\ar[r] & R^0(BGL,1)\ar[r]\ar[d] & R^0(BGL)\ar[r]\ar[d] & R^0\ar[r]\ar[d]     & 0\\
0\ar[r] & \hom_{R_0}(I,R_0)\ar[r]      & \hom_{R_0}(A,R_0)\ar[r]    & \hom_{R_0}(R_0,R_0)\ar[r] & 0}
$$
and
$$
\xymatrix{
0\ar[r] & R^0(BGL^{\times 2},BGL^{\lor 2})\ar[r]\ar[d] & R^0(BGL^{\times 2})\ar[r]\ar[d]          & R^0(BGL^{\lor 2})\ar[r]\ar[d]                & 0\\
0\ar[r] & \hom_{R_0}(I\otimes_{R_0} I,R_0)\ar[r]       & \hom_{R_0}(A\otimes_{R_0} A,R_0)\ar[r]   & \hom_{R_0}(R_0\oplus I^{\oplus 2},R_0)\ar[r] & 0}
$$
of $R^0$-modules.
According to Lemmas \ref{add} and \ref{redadd}, we have an exact sequence
$$
0\lra \Add(BGL,\Omega^\i R)\lra R^0(BGL,1)\lra R^0(BGL^{\times 2},BGL^{\lor 2})
$$
in which the map on the right is the cohomology of the map
$$
\mu-p_1-p_2:(BGL^{\times 2},BGL^{\lor 2})\lra (BGL,1)
$$
($\mu$ is the addition and the $p_i$ are the projections); moreover, this map is the $R_0$-module dual of the multiplication $I\otimes_{R_0}I\to I$.
Hence these short exact sequences assemble into a diagram
$$
\xymatrix{
0\ar[r] & \hom_{R_0}(I/I^2,R_0)\ar[r]\ar[d]     & \Add(BGL,\Omega^\i R)\ar[r]\ar[d] & 0\ar[d]\ar[r]                          & 0\\
0\ar[r] & \hom_{R_0}(I,R_0)\ar[r]\ar[d]         &  R^0(BGL)\ar[r]\ar[d]             & R^0\ar[d]\ar[r]                        & 0\\
0\ar[r] & \hom_{R_0}(I\otimes_{R_0}I,R_0)\ar[r] & R^0(BGL^{\times 2})\ar[r]         & \hom_{R_0}(R_0\oplus I^{\oplus 2},R_0)\ar[r] & 0}
$$
of short exact sequences by the snake lemma.
In particular, we see that $\Add(BGL,\Omega^\i R)$ is naturally identified with the dual $\hom_{R_0}(I/I^2,R_0)$ of the module of indecomposables $I/I^2$.
But $I/I^2\cong M=R_0(BGL_1,1)$, the duality map
$$
R^0(BGL_1,1)\lra \hom_{R_0}(R_0(BGL_1,1),R_0)
$$
is an $R^0$-module isomorphism, and the restriction $R^0(BGL,1)\to R^0(BGL_1,1)$ is dual to the inclusion $M\to I$.
\end{pf}

\section{Algebraic Cobordism}

\subsection{The representing spectrum}
For each natural number $n$, let $V_n\to BGL_n$ denote the universal $n$-plane bundle over $BGL_n$.
Then the Thom spaces
$$
MGL_n:=(V_n,V_n-BGL_n)
$$
come equipped with natural maps
$$
MGL_p\land MGL_q\longrightarrow MGL_{p+q}
$$
defined as the composite of the isomorphism
$$
(V_p,V_p-BGL_p)\land (V_q,V_q-BGL_q)\lra (V_p\times V_q,V_p\times V_q - BGL_p\times BGL_q)
$$
and the map on Thom spaces associated to the inclusion of vector bundles
$$
\xymatrix{
V_p\times V_q\ar[d]\ar[r] & V_{p+q}\ar[d]\\
BGL_p\times BGL_q\ar[r] & BGL_{p+q}.}
$$
Restricting this map of vector bundles along the inclusion $1\times BGL_q\to BGL_p\times BGL_q$ gives a map of Thom spaces
$$
(\AA^p,\AA^p-\AA^0)\land MGL_q\to MGL_{p+q},
$$
and these maps comprise the structure maps of the prespectrum $MGL$.
The associated spectrum is defined by
$$
MGL(p):=\colim_q \Omega^q MGL_{p+q},
$$
as evidently the adjoints
$$
MGL(q)\simeq\colim_r\Omega^r MGL_{q+r}\simeq\colim_r\Omega^{p+r} MGL_{p+q+r}\simeq\Omega^p MGL(p+q)
$$
of the structure maps $\Sigma^p MGL(q)\to MGL(p+q)$ are equivalences.
The last equivalence uses the fact that $\PP^1$ is a compact object of the motivic homotopy category.

\begin{defn}[Voevodsky \cite{Voev98}]
Algebraic cobordism is the motivic cohomology theory represented by the motivic spectrum $MGL$.
\end{defn}

\subsection{Algebraic cobordism is the universal oriented motivic ring spectrum}
Just as in ordinary stable homotopy theory, the Thom classes $\theta_n\in R^n(MGL_n)$ coming from an orientation on a motivic ring spectrum $R$ assemble to give a ring map $\theta:MGL\to R$.
We begin with a brief review of this correspondence.

\begin{prop}[Panin, Pimenov, R\"ondigs \cite{PPR07b}]
Let $R$ be a commutative monoid in the homotopy category of motivic spectra.
Then the set of monoidal maps $MGL\to R$ is naturally isomorphic to the set of orientations on $R$.
\end{prop}

\begin{pf}
The classical analysis of complex orientations on ring spectra $R$ generalizes immediately.
A spectrum map $\theta:MGL\to R$ is determined by a compatible family of maps
$
\theta_n:MGL_n\to R^n,
$
which is to say a family of universal Thom classes $\theta_n\in R^n(MGL_n)$.
An arbitrary $n$-plane bundle $V\to X$, represented by a map $X\to BGL_n$, induces a map of Thom spaces $V/V-X\to MGL_n$, so $\theta_n$ restricts to a Thom class in $R^n(V/V-X)$.
Moreover, these Thom classes are multiplicative and unital precisely when $\theta:MGL\to R$ is monoidal.
Conversely, an orientation on $R$ has, as part of its data, Thom classes $\theta_n\in R^n(MGL_n)$ for the universal bundles $V_n\to BGL_n$, and these assemble to form a ring map $\theta:MGL\to R$.
\end{pf}


Again, just as in topology, an orientation on $R$ is equivalent to a compatible family of $R$-theory Chern classes for vector bundles $V\to X$.
This follows from the Thom isomorphism $R^*(BGL_n)\cong R^*(MGL_n)$.

More difficult is the fact that an orientation on a ring spectrum $R$ is uniquely determined by the first Thom class alone; that is, a class $\theta_1\in R^1(BGL_1)=R$ whose restriction $i^*\theta_1\in R^1(S^1)$ along the inclusion $S^1\to MGL_1$ corresponds to $1\in R^0(S^0)$ via the suspension isomorphism $R^1(S^1)\cong R^0(S^0)$.
This is a result of the splitting principle, which allows one to construct Thom classes (or Chern classes) for general vector bundles by descent from a space over which they split.
See Adams \cite{Ad74} and Panin-Pimenov-R\"ondigs \cite{PPR07b} for details.

\subsection{A ring spectrum equivalent to $PMGL$}
The wedge
$$
\underset{n\in\NN}{\bigvee}\Sigma^\i MGL_n
$$
forms a ring spectrum with unit $\SS\simeq\Sigma^\i MGL_0$ and multiplication
$$
\underset{p}{\bigvee}\Sigma^\i MGL_p\land\underset{q}{\bigvee}\Sigma^\i MGL_q\lra\underset{p,q}{\bigvee}\Sigma^\i MGL_p\land MGL_q\lra\underset{n}{\bigvee}\Sigma^\i MGL_n
$$
induced by the maps $MGL_p\land MGL_q\to MGL_{p+q}$.
Evidently, a (homotopy class of a) ring map $\bigvee_n\Sigma^\i MGL_n\to R$ is equivalent to a family of {\em degree zero} Thom classes $$\vartheta_n\in R^0(MGL_n)$$
with $\vartheta_0=1\in R^0(MGL_0)=R^0$ such that $\vartheta_{p+q}$ restricts via $MGL_p\land MGL_q\to MGL_{p+q}$ to the product $\vartheta_p\vartheta_q$.
This is {\em not} the same as an orientation on $R$, as there is nothing forcing $\vartheta_1\in R^0(MGL_1)$ to restrict to a unit in $R^0(S^1)$.
Clearly we should impose this condition, which amounts to inverting $\beta:\PP^1\to\PP^\i$.

\begin{prop}
A ring map $PMGL\to R$ induces a ring map $\bigvee_n\Sigma^\i MGL_n[1/\beta]\to R$.
\end{prop}

\begin{pf}
A ring map $\theta:PMGL\to R$ consists of a ring map $MGL\to R$ and a unit $\mu\in\pi_1 R$.
This specifies Thom classes $\theta_n\in R^n(MGL_n)$, and therefore Thom classes
$$
\vartheta_n:=\mu^n\theta_n\in R^0(MGL_n)
$$
such that
$$
\vartheta_p\vartheta_q=\mu^{p+q}\theta_p\theta_q=\mu^{p+q}i^*\theta_{p+q}=i^*\vartheta_{p+q}\in R^0(MGL_p\land MGL_q),
$$
where $i$ is the map $MGL_p\land MGL_q\to MGL_{p+q}$.
This gives a ring map $\vartheta:\bigvee_n\Sigma^\i MGL_n\to R$, and therefore the desired map, provided $\beta$ is sent to a unit.
But this is clear:
as a class in $R^0(S^1)$,
$$
\vartheta(\beta)=\beta^*\vartheta_1=\mu\beta^*\theta_1,
$$
and $\beta^*\theta_1\in R^1(S^1)$ is the image of $1\in R^0(S^0)$ under the isomorphism $R^0(S^0)\cong R^1(S^1)$.
\end{pf}

\begin{prop}
The ring map
$
\underset{n\in\NN}{\bigvee}\Sigma^\i MGL_n[1/\beta]\to PMGL
$
is an equivalence.
\end{prop}

\begin{pf}
Write
$$
M:=\underset{n\in\NN}{\bigvee}\Sigma^\i MGL_n[1/\beta]\to PMGL,
$$
and consider the natural transformation of set-valued functors
$$
\rings(PMGL,-)\lra\rings(M,-).
$$
Given a ring spectrum $R$, we have seen that the set $\rings(M,R)$ is naturally isomorphic to the set of collections $\{\vartheta_n\}_{n\in\NN}$ with $\vartheta_n\in R^0(MGL_n)$ such that $\vartheta_{p+q}$ restricts to $\vartheta_p\vartheta_q$, $\vartheta_1$ restricts to a unit in $R^{-1}$, and $\vartheta_0=1\in R^0(S^0)$.
Similarly, the set $\rings(PMGL,R)$ is naturally isomorphic to the product of the set of units in $R^{-1}$ and the set of collections $\{\theta_n\}_{n\in\NN}$ with $R^n(MGL_n)$ such that $\theta_{p+q}$ restricts to $\theta_p\theta_q$, $\theta_1$ restricts to the image of $1\in R^0(S^0)$ in $R^1(S^1)$, and $\theta_0=1\in R^0(S^0)$.

The map $\rings(PMGL,R)\to\rings(M,R)$ sends $\mu\in R^{-1}$ and $\theta_n\in R^n(MGL_n)$ to $\vartheta_n=\mu^n\theta_n$.
We get a natural map back which sends $\vartheta_n\in R^0(MGL_n)$ to $\theta_n=\mu^{-n}\vartheta_n$, where $\mu\in R^{-1}$ in the unit corresponding to $\beta^*\vartheta_1\in R^0(S^1)$.
Clearly the composites are the respective identities, and we conclude that $M\to PMGL$ is an equivalence.
\end{pf}

\subsection{$\Sigma^\i_+ BGL[1/\beta]$ is orientable}
%
Recall from \cite{PPR07b} that, just as in the usual stable homotopy category, an orientation on a ring spectrum $R$ is equivalent to a class in $R^{1}(MGL_1)$ which restricts, under the inclusion $i:S^{1}\to MGL_1$ of the bottom cell, to the class in $R^{1}(S^{1})$ corresponding to the unit $1\in R^{0}(S^{0})$ under the suspension isomorphism $R^{0}(S^{0})\to R^{1}(S^{1})$.
Note also that in the case $R$ is periodic with Bott element $\beta\in R^{0}(S^{1})$, corresponding under the suspension isomorphism to the unit $\mu\in R^{-1}(S^{0})$ with inverse $\mu^{-1}\in R^{1}(S^{0})$, then the suspension isomorphism $R^{0}(S^{0})\to R^{1}(S^{1})$ sends $1$ to $\mu^{-1}\beta$.

Now there's a canonical class $\theta_1\in\Sigma^\i_+ BGL[1/\beta]^{1}(MGL_1)$ such that $$\mu^{-1}\beta=i^*\theta_1\in\Sigma^\i_+ BGL[1/\beta]^{1}(S^{1}).$$
Namely, set
$\theta_1:=\mu^{-1}\vartheta_1$, where $\vartheta_1\in\Sigma^\i_+ BGL[1/\beta]^{0}(MGL_1)$ is the class of the composite
$$
\Sigma^\i MGL_1\simeq\Sigma^\i BGL_1\lra\Sigma^\i_+ BGL\lra\Sigma^\i_+ BGL[1/\beta].
$$
Then $\beta=i^*\mu\theta$, so $\mu^{-1}\beta=i^*\theta$.
%

\begin{prop}\label{can.per.or.}
There is a canonical ring map $\theta:PMGL\to\Sigma^\i_+ BGL[1/\beta]$.
\end{prop}

\begin{pf}
The Thom class $\theta_1\in\Sigma^\i_+ BGL[1/\beta]^0(MGL_1)$ extends, as in \cite{Ad74} or \cite{PPR07b}, to a ring map $MGL\to\Sigma^\i_+ BGL[1/\beta]$, and we have a canonical unit $\mu\in R^{-1}(S^0)$, the image of $\beta\in R^0(S^1)$ under the suspension isomorphism $R^0(S^1)\cong R^{-1}(S^0)$.
\end{pf}

\begin{cor}\label{var.can.per.or}
There is a canonical ring map $\vartheta:\bigvee_n\Sigma^\i MGL_n[1/\beta]\to\Sigma^\i_+ BGL[1/\beta]$.
\end{cor}

\begin{pf}
Precompose the map from the previous Proposition \ref{can.per.or.} with the equivalence
$$\bigvee_n\Sigma^\i MGL_n[1/\beta]\to PMGL.$$
\end{pf}

\subsection{$\vartheta$ is an equivalence}
We analyze the effect of $\vartheta:\bigvee_n\Sigma^\i MGL_n[1/\beta]\to\Sigma^\i_+ BGL[1/\beta]$ on cohomology.
To this end, fix an oriented periodic motivic ring spectrum $R$; we aim to show that the induced map
$$
R^0(\Sigma^\i_+ BGL[1/\beta])\lra R^0(\bigvee_n\Sigma^\i MGL_n[1/\beta])
$$
is an isomorphism.


\begin{lem}\label{algebra}
Let $R$ be a commutative ring and let $A=\colim_{n} A_n$ be a filtered commutative $R$-algebra with the property that $A_p\otimes_R A_q\to A\otimes_R A\to A$ factors through the inclusion $A_{p+q}\to A$ (that is, the multiplication is compatible with the filtration).
Suppose that, for each $n$, the maps $A_{n-1}\to A_n$ are split injections, so that the isomorphisms $A_{n-1}\oplus A_n/A_{n-1}\to A_n$ define an $R$-module isomorphism
$$
\mathrm{gr} A:=\bigoplus_n A_n/A_{n-1}\l{\simeq}{\lra}\colim_n A_n = A
$$
of $A$ with its associated graded.
Then
the multiplication $A_p/A_{p-1}\otimes_R A_q/A_{q-1}\to A_{p+q}/A_{p+q-1}$ makes $\mathrm{gr} A = \bigoplus_n A_n/A_{n-1}$ into a commutative $R$-algebra in such a way that $\mathrm{gr} A\to A$ is an $R$-algebra isomorphism.\hfill$\Box$
\end{lem}

\begin{prop}\label{hopfalgebradecomp}
There is a commuting square of $R^0$-module maps
$$
\xymatrix{
R^0(BGL)\ar[r]\ar[d] & \prod_n R^0(MGL_n)\ar[d]\\
R^0(BGL\times BGL)\ar[r] & \prod_{p,q} R^0(MGL_p\land MGL_q)},
$$
in which the vertical maps are induced by the multiplication on $BGL$ and $\bigvee_n MGL_n$, respectively, and the horizontal maps are isomorphisms.
\end{prop}

\begin{pf}
Set $A:=\colim_n\sym^n_{R_0} R_0(\PP^\i)$, where the map $\sym^{n-1}_{R_0} R_0(\PP^\i)\to\sym^n_{R_0} R_0(\PP^\i)$ is induced by the the inclusion $R_0\cong R_0(\PP^0)\to R_0(\PP^\i)$.
Applying $R_0$ to the cofiber sequence $BGL_{n-1}\to BGL_n\to MGL_n$ yields split short exact sequences
$$
\xymatrix{
\sym^{n-1} R_0(\PP^\i)\ar[r]\ar[d]^\cong & \sym^n R_0(\PP^\i)\ar[r]\ar[d]^\cong & \sym^n R_0(\PP^\i)/\sym^{n-1} R_0(\PP^\i)\ar[d]^\cong\\
R_0(BGL_{n-1})\ar[r]                     & R_0(BGL_n)\ar[r]                     & R_0(MGL_n)}
$$
with $A\cong\colim_n R_0(BGL_n)\cong R_0(BGL)$ a filtered commutative $R$-algebra.
By the lemma, we have a commutative square
$$
\hspace{-.3cm}
\xymatrix{
\underset{p,q}{\bigoplus}\sym^p R_0(\PP^\i)/\sym^{p-1} R_0(\PP^\i)\underset{R_0}{\otimes}\sym^q R_0(\PP^\i)/\sym^{q-1} R_0(\PP^\i)\ar@<4pt>[r]\ar[d] & A\underset{R_0}{\otimes}A\ar[d]\\
\bigoplus_n\sym^n R_0(\PP^\i)/\sym^{n-1} R_0(\PP^\i)\ar[r] & A}
$$
in which the vertical maps are multiplication and the horizontal maps are $R_0$-algebra isomorphisms.
The desired commutative square is obtained by taking $R_0$-module duals.
\end{pf}

\begin{thm}\label{MT0}
The map of oriented periodic motivic ring spectra
$$\vartheta:\bigvee_n\Sigma^\i MGL_n[1/\beta]\to\Sigma^\i_+ BGL[1/\beta]$$
is an equivalence.
\end{thm}

\begin{pf}
We show that the induced natural transformation
$$
\vartheta^*:\rings(\Sigma^\i_+ BGL[1/\beta],-)\lra\rings(\bigvee_n\Sigma^\i MGL_n[1/\beta],-)
$$
is in fact a natural isomorphism.
The result then follows immediately from Yoneda's Lemma.

Fix a ring spectrum $R$, and observe that, for another ring spectrum $A$, $\rings(A,R)$ is the equalizer of the pair of maps from $R^0(A)$ to $R^0(A\land A)\times R^0(\SS)$ which assert the commutativity of the diagrams
$$
\xymatrix{
A\land A\ar[r]\ar[d] & A\ar[d]\\
R\land R\ar[r]       & R}
\hspace{2cm}\text{and}\hspace{2cm}
\xymatrix{
\SS\ar[rd]\ar[r] & A\ar[d]\\
& R}.
$$
Given a map $\beta:\Sigma^1\SS\to A$, the set $\rings(A[1/\beta],R)$ is the equalizer of the pair of maps from $\rings(A,R)\times R^0(\Sigma^{-1}\SS)$ to $R^0(\SS)$ which assert that the ring map $A\to R$ is such that there's a spectrum map $\Sigma^{-1}\SS\to R$ for which the product
$$
\SS\simeq\Sigma^{1}\SS\land \Sigma^{-1}\SS\lra A\land R\lra R\land R\lra R
$$
is equivalent to the unit $\SS\to R$.
Putting these together, we may express $\rings(A[1/\beta],R)$ as the equalizer of natural pair of maps from $R^0(A)\times R^0(\SS^{-1})$ to $R^0(A\land A)\times R^0(\SS)\times R^0(\SS)$.

We therefore get a map of equalizer diagrams
$$
\hspace{-.1cm}
\xymatrix{
R^0(BGL)\times R^0(\Sigma^{-1}\SS)\ar@<2pt>[r]\ar@<-2pt>[r]\ar[d] & R^0(BGL\times BGL)\times R^0(\SS)\times R^0(\SS)\ar[d]\\
\prod_n R^0(MGL_n)\times R^0(\Sigma^{-1}\SS)\ar@<2pt>[r]\ar@<-2pt>[r] & \prod_{p,q} R^0(MGL_p\land MGL_q)\times R^0(\SS)\times R^0(\SS)},
$$
the equalizer of which is $\vartheta^*$.
Now if $R$ does {\em not} admit the structure of a $PMGL$-algebra, then clearly there cannot be any ring maps from either of the $PMGL$-algebras $\bigvee_n\Sigma^\i MGL_n[1/\beta]$ of $\Sigma^\i_+ BGL[1/\beta]$.
Hence we may assume that $R$ is also an oriented periodic ring spectrum, in which case Proposition \ref{hopfalgebradecomp} implies that the vertical maps are isomorphisms.
\end{pf}

\begin{cor}\label{MT1}
The map of periodic oriented motivic ring spectra
$$\theta:PMGL\to\Sigma^\i_+ BGL[1/\beta]$$
is an equivalence.
\end{cor}

\begin{pf}
$\bigvee_n\Sigma^\i MGL_n[1/\beta]\to PMGL$ is an equivalence.
\end{pf}

%

\section{Algebraic $K$-Theory}

\subsection{The representing spectrum}
%
Let $BGL_\ZZ\simeq\ZZ\times BGL$ denote the group completion of the monoid
$$
BGL_\NN:=\coprod_{n\in\NN} BGL_n.
$$
Given a motivic space $X$, write $K^0(X):=\pi_0\map_S(X,BGL_\ZZ)$.
If $S=\spec\ZZ$ and $X$ is a scheme, this agrees with the homotopy algebraic $K$-theory of $X$ as defined by Weibel \cite{Weibel88}, and if in addition $X$ is smooth, this also agrees with Thomason-Trobaugh algebraic $K$-theory of $X$ \cite{TT}; see Proposition 4.3.9 of \cite{MV99} for details.
As the name suggests, homotopy algebraic $K$-theory is a homotopy invariant version of the Thomason-Trobaugh algebraic $K$-theory, and homotopy invariance is of course a prerequisite for any motivic cohomology theory.

It turns out that the motivic space $BGL_\ZZ$, pointed by the inclusion
$$1\simeq BGL_0\lra BGL_\NN\lra BGL_\ZZ,$$
is the zero space of the motivic spectrum $K$ representing (homotopy) algebraic $K$-theory.
This is a direct corollary of the following famous fact.

\begin{prop}[Motivic Bott Periodicity]\label{Bott}
The adjoint
$$
(BGL_\ZZ,BGL_0)\lra\Omega (BGL_\ZZ,BGL_0)
$$
of the map Bott map $\Sigma (BGL_\ZZ,BGL_0)\to (BGL_\ZZ,BGL_0)$ classifying the tensor product of $(L-1)$ and $V$, where $L\to\PP^1$ is the restriction of the universal line bundle and $V\to BGL_\ZZ$ is the universal virtual vector bundle, is an equivalence.
\end{prop}

\begin{pf}
Quillen's projective bundle theorem \cite{Qu73} implies that the tensor product of vector bundles induces an isomorphism
$$
K^0(\PP^1)\underset{K^0}{\otimes} K^0(X)\lra K^0(\PP^1\times X)
$$
of abelian groups.
It follows that there's an isomorphism of split short exact sequences of $K^0(X)$-modules
$$
\xymatrix{
\lambda K^0(X)[\lambda]/(\lambda^2) \ar[r]\ar[d] & K^0(X)[\lambda]/(\lambda^2) \ar[r]\ar[d] & K^0(X)\ar[d]\\
K^0((\PP^1,\PP^0)\land (X,0))\ar[r] & K^0(\PP^1\times X)\ar[r] & K^0(\PP^0\times X)}.
$$
In particular, $K^0((\PP^1,\PP^0)\land (X,0))\cong K^0(X)$ and similarly $K^0((\PP^1,\PP^0)\land (X,1))\cong K^0(X,1)$.
\end{pf}

Define a sequence of pointed spaces $K(n)$ by
$$
K(n):=(BGL_\ZZ,BGL_0)
$$
for all $n\in\NN$.
By Proposition \ref{Bott}, each $K(n)$ comes equipped with an equivalence
$$
K(n)=(BGL_\ZZ,BGL_0)\lra \Omega(BGL_\ZZ,BGL_0)= \Omega K(n+1),
$$
making $K:=(K(0),K(1),\ldots)$ into a motivic spectrum.


\subsection{A map $\Sigma^\i_+ \PP^\i [1/\beta]\to K$}
Let $\beta:\PP^1\to\PP^\i$ be the map classifying the tautological line bundle on $\PP^1$.
We construct a ring map $\Sigma^\i_+ \PP^\i\to K$ which sends $\beta$ to a unit in $K$, thus yielding a ring map $\Sigma^\i_+ \PP^\i [1/\beta]\to K$.
There is a homotopy commutative ring structure on the motivic space $BGL_\ZZ\simeq\Omega^\i K$ in which addition is induced by the sum of vector bundles and multiplication is induced by the tensor product of vector bundles.

Ring maps $\Sigma^\i_+ \PP^\i\to K$ are adjoint to monoidal maps $\PP^\i\to GL_1 K$, the multiplicative monoid of units (up to homotopy) in the ring space $\Omega^\i K\simeq BGL_\ZZ$.
Since $\pi_0 BGL_\ZZ$ contains a copy of $\ZZ$, the multiplicative units contain the subgroup $\{\pm 1\}\to\ZZ$, giving a map
$$\{\pm 1\}\times BGL\lra GL_1 K.$$
But the inclusion $BGL_1\to BGL$ is monoidal with respect to the multiplicative structure on $BGL$, so we get a monoidal map
$$
\PP^\i\simeq BGL_1\lra\{+1\}\times BGL\lra GL_1 K
$$
and therefore a ring map $\Sigma^\i_+ \PP^\i\to K$.

\begin{prop}
The class of the composite
$$
\Sigma^\i S^1\simeq\Sigma^\i(\PP^1,\PP^0)\lra\Sigma^\i_+\PP^\i\lra K
$$
is equal to that of the $K$-theory Bott element $\beta$, i.e. the class of the reduced tautological line bundle $L-1$ on $\PP^1$.
\end{prop}

\begin{pf}
The map $\Sigma^\i_+\PP^\i\to K$ classifies the tautological line bundle on $\PP^\i$, so the pointed version $\Sigma^\i\PP^\i\to K$ corresponds to the reduced tautological line bundle on $\PP^\i$.
This restricts to the reduced tautological line bundle on $\PP^1$.
\end{pf}

\begin{cor}
There's a canonical ring map $\psi:\Sigma^\i_+ \PP^\i [1/\beta]\to K$.\hfill$\Box$
\end{cor}

\subsection{Comparing $R^0(K)$ and $R^0(L)$}
Let $L$ denote the localized motivic ring spectrum
$$
L:=\Sigma^\i_+\PP^\i[1/\beta].
$$
That is, $L$ is the colimit
$$
L=\colim_n\Sigma^{\infty-n}\PP^\i_+
$$
of a telescope of desuspended suspension spectra.
We show that $\psi:L\to K$ is an equivalence by showing the induced map $R^0(L)\to R^0(K)$ is an isomorphism for a sufficiently large class of motivic spectra $R$.

Throughout this section, we will be considering the motivic space $BGL_\ZZ$ (pointed by $\{0\}\times BGL_0$) multiplicatively, as a homotopy commutative monoid with respect to the smash product.
Accordingly, $\Sigma^\i BGL_\ZZ$ is a ring spectrum, and the monoidal map
$$
\PP^\i_+\simeq BGL_0+BGL_1\lra BGL_\NN\lra BGL_\ZZ
$$
gives $\Sigma^\i BGL_\ZZ$ the structure of a homotopy commutative $\Sigma^\i\PP^\i_+$-algebra.
In particular, the Bott element $\beta\in\pi_1\Sigma^\i\PP^\i_+$ determines a Bott element $\beta\in\pi_1\Sigma^\i BGL_\ZZ$ as well as a Bott element $\beta\in\pi_1 K$.

If $R$ is a homotopy commutative ring spectrum equipped with homotopy element $\alpha\in\pi_n R$, we write
$$
\mu(\alpha)=\Sigma^{-n}(\mu\circ(\alpha\land R)):R\lra\Sigma^{-n}R
$$
for the ``multiplication by $\alpha$'' map, the $n$-fold desuspension of the composite
$$
\Sigma^n R\simeq\Sigma^n\SS\land R\l{\alpha\land R}{\lra} R\land R\l{\mu}{\lra} R.
$$
If $\alpha\in\pi_n R$ is a unit, then this map has an inverse $\mu(\alpha)^{-1}:\Sigma^{-n}R\to R$, the $n$-fold desuspension of the multiplication by $\alpha^{-1}$ map $\mu(\alpha^{-1}):R\to\Sigma^n R$.
For our purposes, $R$ will typically admit a periodic orientation, and $\alpha$ will be the image of the Bott element $\beta\in\pi_1 PMGL\cong\pi_1\Sigma^\i_+ BGL[1/\beta]$ under some ring map $PMGL\to R$.

Finally, we also write
$$
\mu(\beta):BGL_\ZZ\lra\Omega BGL_\ZZ
$$
for multiplication by $\beta$ in the homotopy commutative monoid $BGL_\ZZ$ (regarded multiplicatively).
This is $\Omega^\infty$ applied to the multiplication by $\beta$ map $\mu(\beta):K\to\Sigma^{-1}K$ on $K$-theory, and thus it is the equivalence adjoint to the Bott map
$$
\Sigma BGL_\ZZ\lra BGL_\ZZ
$$
The following lemma is formal.

\begin{lem}\label{adlem}
Let $\varepsilon:\Sigma\Omega BGL_\ZZ\to BGL_\ZZ$ denote the counit of the adjunction $(\Sigma,\Omega)$ applied to $BGL_\ZZ$.
Then the composite
$$
\Sigma BGL_\ZZ\l{\Sigma\mu(\beta)}{\lra} \Sigma\Omega BGL_\ZZ\l{\varepsilon}{\lra} BGL_\ZZ
$$
is the Bott map $\Sigma BGL_\ZZ\lra BGL_\ZZ$.
\end{lem}

\begin{pf}
More generally, if $(\Sigma,\Omega)$ is any adjunction and $\beta^*:\Sigma X\to Y$ is a map left adjoint to $\beta_*:X\to\Omega Y$, then $\beta^*=\varepsilon_Y\circ\Sigma\beta_*$.
\end{pf}

\begin{prop}\label{Bottcom}
The square
$$
\xymatrix{
\Sigma^{\infty+1}\PP^\i_+\ar[r]^{\Sigma\mu(\beta)}\ar[d] & \Sigma^\infty\PP^\i_+\ar[d]\\
\Sigma^{\infty+1} BGL_\ZZ\ar[r]^{\Sigma\mu(\beta)}       & \Sigma^\infty BGL_\ZZ},
$$
in which vertical maps come from the inclusion $i:\PP^\i_+\simeq BGL_0 + BGL_1\to BGL_\ZZ$ and the horizontal maps are the Bott maps, commutes up to homotopy.
\end{prop}

\begin{pf}
The inclusion $\Sigma^\i\PP^\i_+\to\Sigma^\i BGL_\ZZ$ is a map of homotopy commutative ring spectra, and the Bott element $\Sigma^\i\PP^1\to\Sigma^\i BGL_\ZZ$ factors through the Bott element $\Sigma^\i\PP^1\to\Sigma^\i\PP^\i_+$.
\end{pf}

\begin{prop}\label{map(K,R)}
Let $R$ be a homotopy commutative $PMGL$-algebra.
Then the space $\map(K,R)$ of maps from $K$ to $R$ is equivalent to the homotopy inverse limit
$$
\map(K,R)\simeq\holim_n\{\cdots\l{f}{\lra}\map(\Sigma^\i BGL_\ZZ,R)\l{f}{\lra}\map(\Sigma^\i BGL_\ZZ,R)\},
$$
where $f=\mu(\alpha)^{-1}\circ\Sigma^{-1}\circ\mu(\beta)$ is the endomorphism of $\map(\Sigma^\i BGL_\ZZ,R)$ which sends a map $x:\Sigma^\i BGL_\ZZ\to R$ to the composite
$$
\Sigma^\i BGL_\ZZ\l{\mu(\beta)}{\lra}\Sigma^{\i-1} BGL_\ZZ\l{\Sigma^{-1} x}{\lra}\Sigma^{-1} R\l{\mu(\alpha)^{-1}}{\lra} R.
$$
\end{prop}

\begin{pf}
In general, for motivic spectra $M$ and $N$,
$$
\map(M,N)\simeq\holim\{\cdots\lra\map(M(1),N(1))\lra\map(M(0),N(0))\},
$$
where the maps send a map $x:M(n)\to N(n)$ to $\Omega x:M(n-1)\simeq\Omega M(n)\to\Omega N(n)\simeq N(n-1)$.
By adjunction, we may rewrite this as
$$
\map(M,N)\simeq\holim\{\cdots\lra\map(\Sigma^\i\Omega^\i\Sigma M,\Sigma N)\lra\map(\Sigma^\i\Omega^\i M,N)\}.
$$
Now $K$ and $R$ are periodic via equivalences $\mu(\beta):K\to\Sigma^{-1} K$ and $\mu(\alpha):R\to\Sigma^{-1} R$, the diagram
$$
\hspace{-.1cm}
\xymatrix{
\map(BGL_\ZZ,\Omega^\i R) \ar[r]^\Omega\ar[d]^\simeq      & \map(\Omega BGL_\ZZ,\Omega^{\i+1} R)\ar[r]^{\mu(\beta)^*}\ar[d]^\simeq & \map(BGL_\ZZ,\Omega^{\i+1} R)\ar[d]^\simeq\\
\map(\Sigma^\i BGL_\ZZ, R)\ar[r]^{\Sigma^\i\varepsilon^*} & \map(\Sigma^{\i+1}\Omega BGL_\ZZ, R)\ar[r]^{\Sigma^{\i+1}\mu(\beta)^*} & \map(\Sigma^{\i+1} BGL_\ZZ, R)},
$$
in which the vertical arrows are adjunction equivalences, commutes, and according to Lemma \ref{adlem} above, the $\Sigma^\i$ applied to the composite $\varepsilon\circ\Sigma\mu(\beta)$ is $\Sigma\mu(\beta):\Sigma^{\i+1} BGL_\ZZ\to\Sigma^\i BGL_\ZZ$.
Hence
$$
\map(K,R)\simeq\holim\{\cdots\lra\map(\Sigma^\i BGL_\ZZ,R)\lra\map(\Sigma^\i BGL_\ZZ,R)\}
$$
is the homotopy inverse limit of the tower determined by the composite
$$
\hspace{-.05cm}
\map(\Sigma^\i BGL_\ZZ,R)\l{\Sigma\mu(\beta)^*}{\lra}\map(\Sigma^{\i+1} BGL_\ZZ, R)\l{\mu(\alpha)^{-1}_*\circ\Sigma^{-1}}{\lra}\map(\Sigma^\i BGL_\ZZ, R).
$$
That is, the endomorphism $f$ of $\map(\Sigma^\i BGL, R)$ sends $x:\Sigma^\i BGL\to R$ to the composite $\mu(\alpha)^{-1}\circ\Sigma^{-1}(x\circ\Sigma\mu(\beta))=\mu(\alpha)^{-1}\circ\Sigma^{-1}(x)\circ\mu(\beta)$, which is to say that $f=\mu(\alpha)^{-1}\circ\Sigma^{-1}\circ\mu(\beta)$.
\end{pf}


\begin{prop}\label{map(L,R)}
Let $R$ be a homotopy commutative $PMGL$-algebra.
Then the space $\map(L,R)$ of maps from $L$ to $R$ is equivalent to the homotopy inverse limit
$$
\map(L,R)\simeq\holim_n\{\cdots\l{g}{\lra}\map(\Sigma^\i\PP^\infty_+,R)\l{g}{\lra}\map(\Sigma^\i\PP^\infty_+,R)\},
$$
where $g=\mu(\alpha)^{-1}\circ\Sigma^{-1}\circ\mu(\beta)$ is the endomorphism of $\map(\Sigma^\i\PP^\i_+,R)$ which sends a map $y:\Sigma^\i\PP^\i_+\to R$ to the composite
$$
\Sigma^\i\PP^\i_+\l{\mu(\beta)}{\lra}\Sigma^{\i-1}\PP^\i_+\l{\Sigma^{-1} y}{\lra}\Sigma^{-1} R\l{\mu(\alpha)^{-1}}{\lra} R.
$$
\end{prop}

\begin{pf}
By definition, $L=\Sigma^\i_+\PP^\i[1/\beta]=\hocolim_n\Sigma^{\infty-n}\PP^\infty_+$, where the map
$$
\Sigma^{\infty-n}\PP^\infty_+\to\Sigma^{\infty-n-1}\PP^\infty_+
$$
is the $n$-fold desuspension of the multiplication by $\beta$ map $\mu(\beta):\Sigma^{\infty}\PP^\infty_+\to\Sigma^{\infty-1}\PP^\infty_+$.
Hence
$$
\hspace{-.15cm}
\map(L,R)\simeq\holim_n\{\cdots\l{\Sigma^{-1}\circ\Sigma\mu(\beta)^*}{\lra}\map(\Sigma^\i\PP^\i_+,\Sigma R)\l{\Sigma^{-1}\circ\Sigma\mu(\beta)^*}{\lra}\map(\Sigma^\i\PP^\i_+,R)\}.
$$
Again, since $R$ is periodic via the multiplication by $\alpha$ map $\mu(\alpha):R\to\Sigma^{-1}R$, we may compose with $\mu(\alpha)^{-1}$ in order to rewrite this as
$$
\map(L,R)\simeq\holim_n\{\cdots\l{g}{\lra}\map(\Sigma^\i\PP^\i_+, R)\l{g}{\lra}\map(\Sigma^\i\PP^\i_+,R)\},
$$
where $g$ is the endomorphism of $\map(\Sigma^\i\PP^\i_+,R)$ which sends the map $y:\Sigma^\i\PP^\i_+\to R$ to the map $g(y)=\mu(\alpha)^{-1}\circ\Sigma^{-1}(y)\circ\mu(\beta)$.
\end{pf}

\begin{cor}\label{commutes}
Let $R$ be a homotopy commutative $PMGL$-algebra.
Then the square
$$
\xymatrix{
\map(\Sigma^\i BGL_\ZZ,R)\ar[r]^f\ar[d]^{i^*} & \map(\Sigma^\i BGL_\ZZ,R)\ar[d]^{i^*}\\
\map(\Sigma^\i\PP^\i_+,R)\ar[r]^g             & \map(\Sigma^\i\PP^\i_+,R)},
$$
in which the vertical map are induced by the inclusion $i:\PP^\i_+\simeq BGL_0+BGL_1\to BGL_\ZZ$, commutes up to homotopy.
In particular, $\psi^*:\map(K,R)\to\map(L,R)$ is the homotopy inverse limit of the map of towers
$$
\hspace{-.15cm}
\xymatrix{
\map(K,R)\ar[r]^\simeq\ar[d]^{\psi^*} & \holim\{\cdots \ar[r]^f & \map(\Sigma^\i BGL_\ZZ, R)\ar[r]^f\ar[d]^{i^*} & \map(\Sigma^\i BGL_\ZZ, R)\}\ar[d]^{i^*}\\
\map(L,R)\ar[r]^\simeq                & \holim\{\cdots \ar[r]^g & \map(\Sigma^\i\PP^\i_+, R)\ar[r]^g             & \map(\Sigma^\i\PP^\i_+, R)\}}
$$
obtained from iterating this commuting square.
\end{cor}

\begin{pf}
This is immediate from Lemmas \ref{map(K,R)} and \ref{map(L,R)}.
\end{pf}

\subsection{A useful splitting}
To complete the analysis of $\psi^*:\map(K,R)\to\map(L,R)$, we must split the space of additive maps from $BGL_\ZZ$ to $\Omega^\i R$ off of the space of all maps from $BGL_\ZZ$ to $\Omega^\i R$.

\begin{prop}\label{loopsadds}
Let $R$ be a motivic spectrum equipped with an equivalence $\mu(\alpha):R\to\Sigma^{-1}R$.
Then given a map $x:\Omega^\i K\to\Omega^\i R$, the map
$$
\Omega^\i\mu(\alpha)^{-1}\circ\Omega(x)\circ\Omega^\i\mu(\beta):\Omega^\i K\lra\Omega^{\i+1} K\lra\Omega^{\i+1} R\lra\Omega^\i R
$$
is a homomorphism for the additive structures on $\Omega^\i K$ and $\Omega^\i R$.
\end{prop}

\begin{pf}
If $X$ is a motivic space equipped with an equivalence $X\to\Omega Y$, then $X$ is a group object in the homotopy category of motivic spaces; if in addition $Y\simeq\Omega Z$, then $X$ is an abelian group object.
In particular, the additions on $\Omega^\i K$ and $\Omega^\i R$ are induced by the equivalences $\Omega^\i\mu(\beta):\Omega^\i K\to\Omega^{\i+1} K$ and $\Omega^\i\mu(\alpha):\Omega^\i R\to\Omega^{\i+1} R$, respectively, and $\Omega^\i\mu(\alpha)^{-1}\circ\Omega(x)\circ\Omega^\i\mu(\beta)$ is a map of loop spaces and therefore respects this addition.
\end{pf}

%

\begin{prop}\label{Split}
Let $R$ be a homotopy commutative $PMGL$-algebra.
Then there exists a canonical section $s:R^{\PP^\i_+}\to R^{BGL_\ZZ}$ of the restriction $r=i^*:R^{BGL_\ZZ}\to R^{\PP^\i_+}$ induced by the inclusion $i:\PP^\i_+\simeq BGL_0+ BGL_1\to BGL_\ZZ$.
\end{prop}

\begin{pf}
Set $Y=R^{BGL_\ZZ}$ and $Z=R^{\PP^\i_+}$.
By Proposition \ref{split}, the additive maps from $BGL_\ZZ$ to $\Omega^\i R$ define a canonical section $\pi_0 Z\to\pi_0 Y$ of the surjection $\pi_0 Y\to\pi_0 Z$.
We must lift this to a map of spectra $s:Z\to Y$.

By Proposition \ref{cohomologyPtensorZ}, we have isomorphisms $R^0(\PP^\i_+\otimes Z)\cong R^0(Z)\otimes_{\pi_0 R}\pi_0 Z\cong Z^0(Z)$.
Combined with the section $\pi_0 Y\to\pi_0 Z$, this induces a map
$$
Z^0(Z)\to R^0(Z)\otimes_{\pi_0 R}\pi_0 Z\to R^0(Z)\otimes_{\pi_0 R}\pi_0 Y\to Y^0(Z).
$$
Take $s\in Y^0(Z)$ to be the image of $1\in Z^0(Z)$ under this map.
\end{pf}

\begin{cor}\label{split}
Let $R$ be a homotopy commutative $PMGL$-algebra.
Then
$$
\map(\Sigma^\i BGL_\ZZ,R)\simeq\map(\Sigma^\i\PP^\i_+,R)\times X
$$
for some space $X$.
\end{cor}

\begin{pf}
Take $X$ to be the global points of the motivic space obtained by applying $\Omega^\i$ to the fiber of $r:R^{BGL_\ZZ}\to R^{\PP^\i_+}$.
\end{pf}

We write
$$
r:\map(\Sigma^\i BGL_\ZZ,R)\lra\map(\Sigma^\i\PP^\i_+,R)
$$
for the restriction and
$$
s:\map(\Sigma^\i\PP^\i_+,R)\lra\map(\Sigma^\i BGL_\ZZ,R)
$$
for the section corresponding to the inclusion of the additive maps from $BGL_\ZZ$ to $\Omega^\i R$ into all maps from $BGL_\ZZ$ to $\Omega^\i R$.
By Corollary \ref{commutes}, we have that $r\circ f = g\circ r$; however, the next proposition shows that in fact $f\simeq s\circ g\circ r$, which is stronger since $r\circ f\simeq r\circ s\circ g\circ r\simeq g\circ r$ as $s$ is a section of $r$.

%

\begin{prop}\label{f=sgr}
Let $R$ be a homotopy commutative $PMGL$-algebra.
Then
$$
f\simeq s\circ g\circ r:\map(\Sigma^\i BGL_\ZZ,R)\lra\map(\Sigma^\i BGL_\ZZ,R).
$$
\end{prop}

\begin{pf}
%
Note that $f$ and $g$ are induced from corresponding maps $f:R^{BGL_\ZZ}\to R^{BGL_\ZZ}$ and $g:R^{\PP^\i_+}\to R^{\PP^\i_+}$, respectively, and that $r$ and $s$ come similarly from maps $r:R^{BGL_\ZZ}\to R^{\PP^\i_+}$ and $s:R^{\PP^\i_+}\to R^{BGL_\ZZ}$.
Thus it's enough to check that
$$
f=s\circ g\circ r\in\pi_0\map(R^{BGL_\ZZ},R^{BGL_\ZZ}).
$$
Since $r\circ f\simeq g\circ r$, we may instead show that $f\simeq s\circ r\circ f$.
By definition,
$$
f\in\pi_0\map(R^{BGL_\ZZ},R^{BGL_\ZZ})\cong R^0(BGL_\ZZ)\underset{R^0}{\widehat{\otimes}} R^0(R^{BGL_\ZZ})
$$
is the image of
$$
1\in\pi_0\map(R^{BGL_\ZZ},R^{BGL_\ZZ})\cong R^0(BGL_\ZZ)\underset{R^0}{\widehat{\otimes}} R^0(R^{BGL_\ZZ})
$$
under the map obtained by applying $(-)\widehat{\otimes}_{R^0} R^0(R^{BGL_\ZZ})$ to the map
$$
\Omega^\i\mu(\alpha)^{-1}_*\Omega^\i\mu(\beta)^*\Omega:\pi_0\map(BGL_\ZZ,\Omega^\i R)\lra\pi_0\map(BGL_\ZZ,\Omega^\i R).
$$
One similarly checks that $r$ and $s$ are obtained by applying $(-)\widehat{\otimes}_{R_0} R^0(R^{BGL_\ZZ})$ to the restriction $R^0(BGL_\ZZ)\to R^0(\PP^\i)$ and its section $R^0(\PP^\i)\to R^0(BGL_\ZZ)$, respectively.
Now, according to Proposition \ref{loopsadds}, as a map of motivic loop spaces, $\Omega^\i\mu(\alpha)^{-1}_*\Omega^\i\mu(\beta)^*\Omega$ sends $x:BGL_\ZZ\to\Omega^\i R$ to the additive map
$$
\Omega^\i\mu(\alpha)^{-1}\circ\Omega(x)\circ\Omega^\i\mu(\beta):BGL_\ZZ\lra\Omega BGL_\ZZ\lra\Omega^{\i+1} R\lra\Omega^\i R,
$$
which is to say that it factors through the inclusion $\Add(BGL_\ZZ,\Omega^\i R)\cong R^0(\PP^\i_+)\to R^0(BGL_\ZZ)$.
Hence $f$ is in the image of
$$
\xymatrix{
R^0(\PP^\i_+)\widehat{\otimes}_{R^0}R^0(R^{BGL_\ZZ})\ar[r]^s\ar[d]^\cong & R^0(BGL_\ZZ)\widehat{\otimes}_{R^0}R^0(R^{BGL_\ZZ})\ar[d]^\cong\\
\pi_0\map(R^{BGL_\ZZ},R^{\PP^\i_+})                 \ar[r]^s             & \pi_0\map(R^{BGL_\ZZ},R^{BGL_\ZZ})},
$$
so $f=s\circ\tilde{f}$ for some $\tilde{f}:R^{BGL_\ZZ}\to R^{\PP^\i_+}$.
But then $s\circ r\circ f\simeq s\circ r\circ s\circ\tilde{f}\simeq s\circ\tilde{f}\simeq f$.
\end{pf}

%

\begin{lem}\label{smalltower}
Suppose given a homotopy commutative diagram of (ordinary) spectra
$$
\xymatrix{
Y\ar[r]^f\ar[d]_r & Y\ar[d]^r\\
Z\ar[r]^g         & Z}
$$
such that $r:Y\to Z$ admits a section $s:Z\to Y$ with
$$f\simeq s\circ g\circ r:Y\to Y.$$
Then the natural map from the homotopy limit of the tower $\{\cdots\to Y\to Y\}$, obtained by iterating $f$, to the homotopy limit of the tower $\{\cdots\to Z\to Z\}$, obtained by iterating $g$, is an equivalence.
\end{lem}

\begin{pf}
Since $Z$ is a retract of $Y$, we may write $Y\simeq Z\times X$ for some spectrum $X$ such that the fiber of $f$ over $g$ is the trivial map $X\to X$.
Now consider the diagram
$$
\xymatrix{
W                         \ar[r]\ar[d] & \prod_n X\ar[r]      \ar[d] & \prod_n X\ar[d]\\
\holim\{\cdots\to Y\to Y\}\ar[r]\ar[d] & \prod_n Y\ar[r]^{1-f}\ar[d] & \prod_n Y\ar[d]\\
\holim\{\cdots\to Z\to Z\}\ar[r]       & \prod_n Z\ar[r]^{1-g}       & \prod_n Z}
$$
in which the rows and columns are fiber sequences.
Then the fiber of $1-f$ over $1-g$ is the identity $\prod_n X\to\prod_n X$, so $W$ is trivial and $\holim\{\cdots\to Y\to Y\}\simeq\holim\{\cdots\to Z\to Z\}$.
\end{pf}

\begin{prop}\label{R0K}
Let $R$ be an orientable motivic ring spectrum and let $\alpha\in\pi_1 R$ be a unit.
Then
$$
\psi^*:\map(K,R)\lra\map(L,R)
$$
is an equivalence.
\end{prop}

\begin{pf}
Set $Y=\map(\Sigma^\i BGL_\ZZ,R)$ and $Z=\map(\Sigma^\i\PP^\i_+,R)$, so that $Y\simeq Z\times X$ via $r:Y\to Z$ and its section $s:Z\to Y$.
By Proposition \ref{f=sgr}, $f\simeq s\circ g\circ r:Y\to Y$, and since $r$ and $s$ are infinite loop maps we may regard them as maps of (ordinary) connective spectra.
The result then follows from Lemma \ref{smalltower}.
\end{pf}

\begin{cor}\label{R(K)=R(L)}
Then $\psi:L\to K$ induces an isomorphism $\psi^*:R^0(K)\to R^0(L)$ for any orientable periodic motivic ring spectrum $R$.
\end{cor}

\begin{pf}
This is immediate from Proposition \ref{R0K} above, since the spectrum of motivic spectrum maps from $L$ to $R$ admits precisely the same description as that of the spectrum of motivic spectrum maps from $K$ to $R$.
Indeed,
$$
L\simeq\colim_n\{\Sigma^\i_+\PP^\i\l{\beta_*}{\lra}\Sigma^{-1}\Sigma^\i_+\PP^\i\l{\beta_*}{\lra}\cdots\},
$$
and we see that
$$
R^L\simeq\holim_n\{\cdots\l{g}{\lra} R^{\PP^\i_+}\l{g}{\lra} R^{\PP^\i_+}\},
$$
where the map $g:R^{\PP^\i_+}\to R^{\PP^\i_+}$ sends $\lambda:\Sigma^\i_+\PP^\i\to R$ to $\alpha_*^{-1}\circ\Sigma^{-1}\lambda\circ\beta_*$, the composite
$$
\Sigma^{\i}_+\PP^\i\to\Sigma^{-1}\Sigma^\i_+\PP^\i\to\Sigma^{-1}R\to R,
$$
just as above.
\end{pf}

By the homotopy category of orientable periodic spectra, we mean the full subcategory of the homotopy category of spectra on the orientable and periodic objects.
In other words, $R$ is an orientable periodic spectrum if there exists a homotopy commutative ring structure on $R$ which admits a ring map $PMGL\to R$.
Note that, according to this definition, maps between orientable periodic spectra need not preserve potential orientations or even ring structures.

\begin{prop}\label{natiso}
$\psi$ induces an isomorphism $\psi^*:[K,-]\to [L,-]$ of functors from the homotopy category of orientable periodic spectra to abelian groups.
\end{prop}

\begin{pf}
Let $R$ be an orientable periodic spectrum.  Then
$$
\psi^*:[K,R]=R^0(K)\to R^0(L)=[L,R]
$$
is an isomorphism by Corollary \ref{R(K)=R(L)}, and this isomorphism is natural in spectrum maps $R\to R'$, provided of course that $R'$ is also orientable and periodic.
\end{pf}

\begin{thm}\label{mt2}
The ring map $\psi:L\to K$ is an equivalence.
\end{thm}

\begin{pf}
Let $\varphi^*:[L,-]\to [K,-]$ be the inverse of the isomorphism $\psi^*:[K,-]\to [L,-]$ of Proposition \ref{natiso}, and let $\varphi:K\to L$ be the map obtained by applying $\varphi^*$ to the identity $1\in [L,L]$.
It follows from the Yoneda lemma that $\varphi^*$ is precomposition with $\varphi$.
The equations $\varphi^*\circ\psi^*=1_K^*$ and $\psi^*\circ\varphi^*=1_L^*$ imply that $\psi\circ\varphi=1_K$ and $\varphi\circ\psi=1_L$ in the homotopy category of orientable periodic spectra, and therefore that $\psi\circ\varphi=1_K$ and $\varphi\circ\psi=1_L$ in the homotopy category of spectra.
Hence $\psi:L\to K$ is an equivalence with inverse $\varphi:K\to L$.
\end{pf}

\section{Applications}

\subsection{The motivic Conner-Floyd theorem}
The classical theorem of Conner and Floyd shows that complex cobordism determines complex $K$-theory by base change.
More precisely, writing $PMU$ for periodic complex cobordism and $KU$ for complex $K$-theory, then, for any finite spectrum $X$, the natural map
$$
PMU^0(X)\otimes_{PMU^0}KU^0\to KU^0(X)
$$
is an isomorphism of $KU^0$-modules.

\begin{re}
This is the precursor of the more general notion of Landweber exactness.
In \cite{La}, P. Landweber gives a necessary and sufficient condition on an $MU_*$-module $G$ so that the functor
$
(-)\otimes_{MU_*} G,
$
from $(MU_*,MU_*MU)$-comodules to graded abelian groups, is exact.
For $G=K_*$, it follows that the natural map
$$
MU^*(-)\otimes_{MU^*}K^*\to K^*(-)
$$
is an isomorphism.
\end{re}

%
%
%

We now turn to the motivic version of the theorem of Conner and Floyd.
A motivic spectrum $X$ is said to be {\em compact} if $[X,-]$, viewed as a functor from motivic spectra to abelian groups, commutes with filtered colimits.

%

\begin{prop}\label{cfepi}
Let $X$ be a compact motivic spectrum.
Then the natural map
$$
PMGL^0(X)\underset{PMGL^0}{\otimes}K^0\lra K^0(X)
$$
is surjective.
\end{prop}

\begin{pf}
Set $B:=\Sigma^\i_+ BGL$ and $A:=\Sigma^\i_+\PP^\i$.
Then the determinant map $r:B\to A$ admits a section $s:A\to B$, so, for each $n$, $\Sigma^{-n}A$ is a retract of $\Sigma^{-n}B$ and $\Sigma^{-n}B^0(X)\to\Sigma^{-n}A^0(X)$ is surjective.
Since $X$ is compact, the colimit
$$
B[1/\beta]^0(X)\cong\colim_n \Sigma^{-n}B^0(X)\longrightarrow\colim_n \Sigma^{-n}A^0(X)\cong A[1/\beta]^0(X)
$$
is also surjective, and we see that
$$
B[1/\beta]^0(X)\underset{B[1/\beta]^0}{\otimes} A[1/\beta]^0\longrightarrow A[1/\beta]^0(X)\underset{B[1/\beta]^0}{\otimes} A[1/\beta]^0\cong A[1/\beta]^0(X)
$$
is surjective as well.
\end{pf}

\begin{thm}
Let $X$ be a compact motivic spectrum.
Then the natural map
$$
PMGL^{0}(X)\underset{PMGL^{0}}{\otimes} K^{0}\lra K^{0}(X)
$$
is an isomorphism.
\end{thm}

\begin{pf}
According to Proposition \ref{cfepi}, it's enough to show that the map is injective.
For simplicity of notation, set $R:=PMGL$, define a contravariant functor $J^0(-)$ from compact motivic spectra to $R^0$-modules by the rule
$$
J^0(X):=\ker\{R^0(X)\to K^0(X)\},
$$
and write $J^0$ for $J^0(\SS)$.
Since the tensor product is right exact, the map
$$
J^0(X)\otimes_{R^0}K^0\lra\ker\{R^0(X)\otimes_{R^0}K^0\to K^0(X)\otimes_{R^0}K^0\}
$$
is surjective, so in light of the isomorphism
$$
K^0(X)\otimes_{R^0}K^0\cong K^0(X)\otimes_{K^0}K^0\cong K^0(X)
$$
it's enough to show that $J^0(X)\otimes_{R^0}K^0$ is zero, or, equivalently, that
$$
J^0(X)\otimes_{R^0}J^0\lra J^0(X)\otimes_{R^0}R^0\cong J^0(X)
$$
is surjective.
To this end, set
$$
I^0(X):=\im\{J^0(X)\otimes_{R^0}J^0\to J^0(X)\otimes_{R^0}R^0\cong J^0(X)\};
$$
we must show that $I^0(X)\cong J^0(X)$.

Now, writing $B:=\Sigma^\i_+ BGL$ and $A:=\Sigma^\i_+\PP^\i$ as above, and using the compactness of $X$, we see that any element of
\begin{align*}
J^0(X)&\cong\ker\{\colim_n[X,\Sigma^{-n}B]\to\colim_n[X,\Sigma^{-n}A]\}\\
      &\cong\colim_n\ker\{[X,\Sigma^{-n}B]\to[X,\Sigma^{-n}A]\}
\end{align*}
is represented by a map
$$
x:\Sigma^n X\to B\simeq\colim_p\colim_q\Sigma^\i_+\Grass_{p,q},
$$
which, by compactness, factors as $f_x:\Sigma^n X\to Y_x$ followed by $y:Y_x\to B$ for $Y_x\simeq\Sigma^\i_+\Grass_{p,q}$ the suspension spectrum of a finite Grassmannian.
This yields a commuting diagram
$$
\xymatrix{
\Sigma^n X\ar[r]^x\ar[d]_{f_x} & B\ar[d]^{r}\\
Y_x\ar[r]\ar[ur]^{y}           & A}
$$
in which $r\circ x$ is trivial and the determinant map $r:B\to A$ admits a section $s:A\to B$.
Of course, as $r\circ y$ need not be trivial, set $y':=y-s\circ r\circ y$, so that
$$
y'\circ f\simeq (y-s\circ r\circ y)\circ f\simeq y\circ f-s\circ r\circ y\circ f\simeq x-s\circ r\circ x\simeq x
$$
and $r\circ y'\simeq 0$, which is to say that $y'\in J^0(X)$ and $f_x^*y'=x$.

Finally, according to Proposition \ref{cograss}, $R^0(Y_x)\otimes_{R^0}K^0\cong K^0(Y_x)$ for each $x\in J^0(X)$, so we must have surjections
$$
J^0(Y_x)\otimes_{R^0} J^0\lra J^0(Y_x)\otimes_{R^0} R^0\cong J^0(Y_x).
$$
Adding these together, we obtain a morphism of short exact sequences
$$
\xymatrix{
0\ar[r] & \bigoplus_x I^0(Y_x)\ar[r]\ar[d] & {\bigoplus}_x J^0(Y_x)\ar[r]\ar[d]^{\bigoplus_x f_x^*} & 0 \ar[r]\ar[d]                & 0\\
0\ar[r] & I^0(X)              \ar[r]       & J^0(X)                \ar[r]                           & J^0(X)\otimes_{R^0} K^0\ar[r] & 0}
$$
such that $\bigoplus_x f_x^*:\bigoplus_x J^0(Y_x)\to J^0(X)$ is surjective.
It follows that $I^0(X)\cong J^0(X)$.
\end{pf}

\subsection{$PMGL$ and $K$ are $E_\i$ motivic spectra}
As a final application, we show that $PMGL$ and $K$ are $E_\i$ motivic spectra.
As we shall see, this is an immediate consequence of the fact that $PMGL$ and $K$ are obtained through a localization of the category of $E_\i$ $R$-algebras for some $E_\i$ motivic spectrum $R$.
Roughly, given an element $\beta\in\pi_{p,q} R$, the functor which sends the $R$-module $M$ to $M[1/\beta]:=M\land_R R[1/\beta]$ defines a monoidal localization of the category of $R$-modules, so it extends to a localization of the category of $E_\i$ $R$-algebras.
Taking $R=\Sigma^\i_+ BGL$ and $\beta$ the Bott element, we see that $PMGL$ is the localization of the initial $E_\i$ $R$-algebra and $K$ is the localization of the determinant $E_\i$ $R$-algebra.

In order to make this precise, we fix a suitable symmetric monoidal model category $(\mod_\SS,\land_\SS)$ of motivic spectra, such as motivic $\SS$-modules \cite{Hu} or motivic symmetric spectra \cite{Jard00}.
For sake of definiteness, we adopt the formalism of the latter; nevertheless, we refer to motivic symmetric spectra as $\SS$-modules, as they are indeed modules over the symmetric motivic sphere $\SS$.

Recall that a {\em motivic symmetric sequence} is a functor from the groupoid $\Sigma$ of finite sets and isomorphisms to {\em pointed} motivic spaces.
It is sometimes convenient to use a skeleton of $\Sigma$, so we simply write $n$ for a finite set with $n$ elements and $\Sigma(n)$ for its automorphism group.
Motivic symmetric sequences form a symmetric monoidal category under the smash product defined by
$$
(X\land Y)(n):=\bigvee_{n=p+q}\Sigma(n)_+\land_{\Sigma(p)\times\Sigma(q)} X(p)\land Y(q).
$$
The motivic sphere $\SS$ has a natural interpretation as the motivic symmetric sequence in which $\SS(n)$ is the pointed $\Sigma(n)$-space associated to the pair $(\AA^n,\AA^n-\AA^0)$, where $\Sigma(n)$ acts by permutation of coordinates.
The $\Sigma(p)\times\Sigma(q)$-equivariant maps
$$
(\AA^p,\AA^p-\AA^0)\land(\AA^q,\AA^q-\AA^0)\lra(\AA^{p+q},\AA^{p+q}-\AA^0)
$$
give $\SS$ the structure of a commutative monoid for this smash product.
An {\em $\SS$-module} is then a motivic symmetric sequence equipped with an action of $\SS$, which is to say a sequence $X(p)$ of pointed $\Sigma(p)$-equivariant motivic spaces equipped with $\Sigma(p)\times\Sigma(q)$-equivariant maps
$$
(\AA^p,\AA^p-\AA^0)\land X(q)\lra X(p+q);
$$
the fact that $\SS$ is a commutative monoid implies that $\land$ extends to a smash product $\land_\SS$ on the category $\mod_\SS$ of $\SS$-modules.
There are monoidal functors
$$
\hspace{-.15cm}
\{\text{Motivic spaces}\}\!\to\!\{\text{Motivic symmetric spaces}\}\!\to\!\{\text{Motivic symmetric spectra}\}
$$
in which the righthand map is the free $\SS$-module functor, the left hand map sends the motivic space $X$ to the constant motivic symmetric space $X_+$, and the composite is a structured version of the suspension spectrum functor $\Sigma^\i_+$.

\begin{prop}\label{strictcom}
The $\SS$-modules $\Sigma^\i_+ BGL$ and $\Sigma^\i_+\PP^\i$ are equivalent to strictly commutative $\SS$-algebras in such a way that the determinant map $\Sigma^\i_+ BGL\to\Sigma^\i_+\PP^\i$ is equivalent to a map of strictly commutative $\SS$-algebras.
\end{prop}

\begin{pf}
For each $n$, write $GL(n)$ for the group $S$-scheme of linear automorphisms of $\AA^n$.
Then $\Sigma(n)$ acts on $GL(n)$ by conjugation via the embedding $\Sigma(n)\to GL(n)$, so that $GL(n)$ is the value at $n$ of a symmetric sequence $GL$ in group $S$-schemes such that the determinant map $GL\to GL_1$ is a morphism of symmetric sequences in group $S$-schemes, where we regard $GL_1$ as a constant symmetric sequence.
Taking classifying spaces, we obtain a morphism of commutative monoid symmetric sequences $BGL\to BGL_1$ in unpointed motivic spaces.
Now let $\SS[BGL]$ and $\SS[BGL_1]$ denote the $\SS$-modules defined by
$$\SS[BGL](n):=\SS(n)\land BGL(n)_+\qquad\text{and}\qquad\SS[BGL_1](n):=\SS(n)\land BGL_{1+},$$
respectively, where $\Sigma(n)$ acts diagonally; note that $\SS[BGL_1]$ is the free $\SS$-module on the motivic symmetric sequence $BGL_{1+}$, whereas the action of $\SS$ on $\SS[BGL]$ is induced by the canonical $\Sigma(p)\times\Sigma(q)$-equivariant inclusions
$$
BGL(q)\lra BGL(p)\times BGL(q)\to BGL(p+q)
$$
coming from the fact that $BGL(p)$ has a canonical basepoint which is fixed by the action of $\Sigma(p)$.
The monoidal structure on $BGL_1$ extends to a strictly commutative $\SS$-algebra structure on $\SS[BGL_1]$, and the strictly commutative $\SS$-algebra structure on $\SS[BGL]$ comes from $\Sigma(p)\times\Sigma(q)$-equivariant pairing
$$
\SS(p)\land BGL(p)_+\land\SS(q)\land BGL(q)_+\lra\SS(p+q)\land BGL(p+q)_+;
$$
moreover, it is clear that the determinant map $\SS[BGL]\to\SS[BGL_1]$ is monoidal with respect to these multiplicative structures.
Hence we are done, provided the underlying ordinary motivic spectra (obtained by forgetting the actions of the symmetric groups) of $\SS[BGL]$ and $\SS[BGL_1]$ are equivalent to the motivic spectra $\Sigma^\i_+ BGL$ and $\Sigma^\i_+\PP^\i$, respectively.
This is immediate for $\SS[BGL_1]$, whose underlying spectrum is the suspension spectrum $\Sigma^\i_+ BGL_1$; on the other hand, the underlying spectrum of $\SS[BGL]$ is the prespectrum $\{S^n\land BGL_{n+}\}$.
But the motivic spectrum associated to $\Sigma^\i_+ BGL$ is given by
$$
\colim_p\Omega^p\Sigma^p\colim_q BGL_q\simeq\colim_p\colim_q\Omega^p\Sigma^p BGL_q\simeq\colim_n\Omega^n\Sigma^n BGL_n,
$$
so the two motivic prespectra are stably equivalent.
\end{pf}

Instead of considering localization in the context of symmetric monoidal model categories, it will be enough to consider localization in the context of symmetric monoidal $\i$-categories, in the sense of Lurie \cite{DAGIII}.
Indeed, if $(\M,\otimes)$ is a symmetric monoidal model category, then, in the notation of \cite{DAGIII}, the commutative algebra objects of the associated symmetric monoidal $\i$-category $\mathrm{N}(\M,\otimes)^\circ$ correspond to {\em coherently homotopy commutative}, or $E_\i$, objects of $(\M,\otimes)$.
Roughly, this is because the free $E_\i$ algebra monad on $(\M,\otimes)$ is a model for the free commutative algebra monad on $\mathrm{N}(\M,\otimes)^\circ$, and consequently $E_\i$ algebras in $(\M,\otimes)$ model commutative algebras in $\mathrm{N}(\M,\otimes)$.
Here $\mathrm{N}$ denotes the simplicial nerve of a simplicial category; if the simplicial category comes from a symmetric monoidal model category, then its simplicial nerve is symmetric monoidal as an $\i$-category.
See \cite{HTT} for facts about $\i$-categories and simplicial nerves, and \cite{DAGIII} for a treatment of commutative algebra in the $\i$-categorical context.

Recall (cf. \cite{HTT}) that a map $F:\C\to\D$ of $\i$-categories is said to be a {\em localization} if $F$ admits a fully faithful right adjoint $G$.
In this case, it is common to identify $\D$ with the full subcategory of $\C$ consisting of those objects in the essential image of $G$ (the ``local objects''), and suppress $\D$ and $G$ from the notation by writing $L$ for the composite $G\circ F:\C\to\D\to\C$.
If $\C$ is the underlying $\i$-category of a symmetric monoidal $\i$-category $(\C,\otimes)$, then we may ask when a localization $L:\C\to\C$ extends to a lax symmetric monoidal functor on $(\C,\otimes)$.
Since maps $L(-)\otimes L(-)\to L((-)\otimes (-))$ correspond to maps $L(L(-)\otimes L(-))\to L((-)\otimes (-))$, and there is a canonical map $L((-)\otimes (-))\to L(L(-)\otimes L(-))$ in the other direction, it is enough to require that this is an equivalence.
Or we may take advantage of the symmetry of the monoidal structure and simply require that
$
L((-)\otimes (-))\lra L((-)\otimes L(-))
$
is an equivalence.

\begin{defn}[\cite{DAGIII}, Definition 1.28]\label{Lcompatible}
Let $(\C,\otimes)$ be a symmetric monoidal $\i$-category and let $L:\C\to\C$ be a localization of the underlying $\i$-category.
Then $L$ is said to be {\em compatible} with $\otimes$ if, for all objects $A$ and $B$ of $\C$, the map
$$
L(A\otimes B)\lra L(A\otimes LB)
$$
is an equivalence.
\end{defn}

\begin{prop}[\cite{DAGIII}, Proposition 1.31]\label{Llaxmonoidal}
Let $(\C,\otimes)$ be a symmetric monoidal $\i$-category, let $L:\C\to\C$ be a localization of the underlying $\i$-category, and suppose that $L$ is compatible with $\otimes$.
Then $L$ extends to a lax symmetric monoidal functor
$$
(L,\otimes):(\C,\otimes)\to(\C,\otimes).
$$
In particular, $L$ preserves algebra and commutative algebra objects of $(\C,\otimes)$.
\end{prop}

Let $\mathrm{N}(\mod_\SS,\land_\SS)^\circ$ denote the symmetric monoidal $\i$-category which arises as the simplicial nerve of the symmetric monoidal simplicial model category $(\mod_\SS,\land_\SS)$ of $\SS$-modules.
Since commutative algebra objects of $\mathrm{N}(\mod_\SS,\land_\SS)^\circ$ are modeled by algebras over a suitable $E_\i$ operad, we refer to commutative algebra objects of $\mathrm{N}(\mod_\SS,\land_\SS)^\circ$ as {\em $E_\i$ $\SS$-algebras}.
Given an $E_\i$ $\SS$-algebra $R$, we write $(\mod_R,\land_R)$ for the resulting symmetric monoidal $\i$-category of $R$-modules, and refer to commutative algebra objects of $(\mod_R,\land_R)$ as {\em $E_\i$ $R$-algebras}.

\begin{prop}\label{smashing}
Let $R$ be an $E_\i$ $\SS$-algebra, let $f\in\pi_{p,q} R$ be an arbitrary element, and write $L_f:\mod_R\to\mod_R$ for the functor which sends the $R$-module $M$ to the $R$-module
$$
M[1/f]:=M\land_R R[1/f].
$$
Then $L_f$ is a localization functor which is compatible with the symmetric monoidal structure $\land_R$ on $\mod_R$;
in particular, $L_f$ extends to a lax monoidal functor $L_f:(\mod_R,\land_R)\to(\mod_R,\land_R)$.
\end{prop}

\begin{pf}
Say that an $R$-module $M$ is $f$-local if the multiplication by $f$ map $f_*:M\to\Sigma^{-p,-q}M$ is an equivalence.
Given an $f$-local $R$-module $M$, the induced map
$$\Map(R[1/f],M)\simeq\lim\{M\leftarrow\Sigma^{p,q}M\leftarrow\cdots\}\simeq M$$
is an equivalence, so
$$
\map(N[1/f],M)\simeq\map(N,\Map(R[1/f],M))\simeq\map(N,M)
$$
is an equivalence for any $R$-module $N$.
Hence $L_f$ is left adjoint to the inclusion of the full subcategory of $f$-local $R$-modules, and is therefore a localization.
Moreover, it is compatible with $\land_R$, since we have equivalences
$$
\hspace{-.15cm}
L_f(M\land_R N)\simeq M\land_R N\land_R R[1/f]\simeq M\land_R N\land_R R[1/f]\land_R R[1/f]\simeq L_f(M\land L_f N).
$$
Hence, by Proposition \ref{Llaxmonoidal}, $L_f$ extends to a lax symmetric monoidal endofunctor $(L_f,\land_R)$ of $(\mod_R,\land_R)$.
\end{pf}

\begin{cor}\label{loc}
Let $R$ be an $E_\i$ $\SS$-algebra and let $f\in\pi_{p,q}$ be a fixed element.
Then $R[1/f]$ is an $E_\i$ $R$-algebra, and therefore also an $E_\i$ $\SS$-algebra.
\end{cor}

\begin{pf}
By Proposition \ref{smashing}, $L_f$ is a lax symmetric monoidal functor with $L_f R\simeq R[1/f]$.
Since lax symmetric monoidal functors preserve commutative algebra objects, we see that $R[1/f]$ is a commutative algebra object in $(\mod_R,\land_R)$.
Lastly, as the forgetful functor from $R$-modules to $\SS$-modules is lax symmetric monoidal, it follows that $R[1/f]$ is also an $E_\i$ $\SS$-algebra.
\end{pf}

\begin{cor}
$MGL$, $PMGL$ and $K$ are $E_\i$ $\SS$-algebras.
\end{cor}

\begin{pf}
The $MGL$ case is already well-known (cf. \cite{Hu}, for example).
By Proposition \ref{strictcom}, $\Sigma^\i_+ BGL$ and $\Sigma^\i_+\PP^\i$ are equivalent to strictly commutative $\SS$-algebras, so they are naturally commutative algebra objects in the symmetric monoidal $\i$-category $\mathrm{N}(\mod_\SS,\land_\SS)$.
Applying Proposition \ref{loc}, we see that $PMGL\simeq\Sigma^\i_+ BGL[1/\beta]$ is an $E_\i$ $\Sigma^\i_+ BGL$-algebra, and likewise that $K\simeq\Sigma^\i_+\PP^\i[1/\beta]$ is an $E_\i$ $\Sigma^\i_+\PP^\i$-algebra.
In particular, $PMGL$ and $K$ are $E_\i$ $\SS$-algebras.
\end{pf}

\begin{prop}
$K$ is an $E_\i$ $PMGL$-algebra.
\end{prop}

\begin{pf}
Note that the Bott element $\Sigma\SS\to\Sigma^\i_+\PP^\i$ factors as the composite of the Bott element $\Sigma\SS\to\Sigma^\i_+ BGL$ followed by the determinant map $\Sigma^\i_+ BGL\to\Sigma^\i_+\PP^\i$.
By Proposition \ref{strictcom}, the determinant map $\Sigma^\i_+ BGL\to\Sigma^\i_+\PP^\i$ is a map of $E_\i$ $\Sigma^\i_+ BGL$-algebras, so by Proposition \ref{loc}, the localization
$
K\simeq\Sigma^\i_+\PP^\i[1/\beta]
$
is an $E_\i$ algebra over $PMGL\simeq\Sigma^\i_+ BGL[1/\beta]$.
%
\end{pf}

\newpage




\begin{thebibliography}{99}


\bibitem{Ad74} J.F. Adams. {\em Stable Homotopy and Generalised Homology}. University of Chicago Press (1974).

\bibitem{Ad76b} J.F. Adams. Primitive elements in the  K-theory of $BSU$. Quart. J. Math. Oxford (2) 27 (1976) no. 106, 253-262.

\bibitem{Arth83} R.D. Arthan. Localization of stable homotopy rings. Math. Proc. Camb. Phil. Soc. 93 (1983) 295-302.

\bibitem{At68} M.F. Atiyah.  {\em K-Theory}. Benjamin (1968).

\bibitem{DI} D. Dugger and D. Isaksen. Motivic cell structures. Algebraic and Geometric Topology, Vol. 5 (2005) 615--652.

\bibitem{Hu} P. Hu.  $S$-modules in the category of schemes.  Mem. Amer. Math. Soc. 161 (2003), no. 767.

\bibitem{Jard00} J.F. Jardine. Motivic symmetric spectra. Documenta Math. 5 (2000) 445-552.

\bibitem{La} P. Landweber. Homological properties of comodules over MU*(MU) and BP*(BP). American Journal of Mathematics, Vol. 98, No. 3 (1976) 591-610.

\bibitem{Le04} M. Levine. A survey of algebraic cobordism. Proc. International Conf. on Algebra; Algebra Colloq. 11 (2004) no.1 79-90.

\bibitem{Le} M. Levine. Algebraic cobordism II.  Preprint.  http://www.math.uiuc.edu/K-theory/0577

\bibitem{LeMo01} M. Levine and F. Morel. Cobordisme alg\'{e}brique I, II. C.R. Acad. Sci. Paris S\'{e}r. I Math.  (8) 332 (2001) 723-728 and (9) 332 (2001) 815-820.

\bibitem{LeMo02} M. Levine and F. Morel. Algebraic cobordism I.  Preprint. http://www.math.uiuc.edu/K-theory/0547

\bibitem{Lur} J. Lurie. Survey article on elliptic cohomology. Preprint. http://www-math.mit.edu/ $^\sim$lurie/papers/survey.pdf

\bibitem{DAGI} J. Lurie. Derived algebraic geometry I: stable $\i$-categories. Preprint. http://arxiv.org/abs/math/0608228

\bibitem{DAGII} J. Lurie. Derived algebraic geometry II: noncommutative algebra. Preprint. http://arxiv.org/abs/math/0702299

\bibitem{DAGIII} J. Lurie. Derived algebraic geometry III: commutative algebra. Preprint. http://arxiv.org/abs/math/0703204

\bibitem{HTT} J. Lurie.  Higher topos theory. Preprint.  http://arxiv.org/abs/math/0608040

\bibitem{Mor04}   F. Morel.  On the motivic $\pi_{0}$ of the sphere spectrum. {\em Axiomatic, Enriched and Motivic Homotopy Theory}
(ed. J.P.C. Greenlees) NATO Science Series II \#131 (2004) 219-260

\bibitem{Mor06}   F. Morel.  {\em Homotopy Theory of Schemes}, (trans. J.D. Lewis). A.M. Soc. Texts and Monographs \#12 (2006).

\bibitem{MV99} F. Morel and V. Voevodsky. ${\mathbb A}^{1}$ homotopy theory of schemes.  Publ. IHES 90 (1999) 45-143.

\bibitem{NOS} N. Naumann, M. Spitzweck and P.A. {\O}stv{\ae}r.  Motivic Landweber exactness. Preprint.  http://arxiv.org/abs/0806.0274

\bibitem{PPR07}  I. Panin, K. Pimenov and O. R\"{o}ndigs. On Voevodsky's algebraic K-theory spectrum $BGL$. Preprint. http://www.math.uiuc.edu/K-theory/0838

\bibitem{PPR07b}  I. Panin, K. Pimenov and O. R\"{o}ndigs. A universality theorem for Voevodsky's algebraic cobordism spectrum. Homology, Homotopy and Applications, Vol. 10 (2008), No.2, pp. 211-226.

\bibitem{PPRcf}  I. Panin, K. Pimenov and O. R\"{o}ndigs.  On the relation of Voevodsky's algebraic cobordism to Quillen's K-theory.  Invent. Math. 175 (2009), No. 2, 435-451.

Preprint. http://www.math.uiuc.edu/K-theory/0847

\bibitem{Qu69} D.G. Quillen. On the formal group laws of unoriented and complex cobordism theory. Bull. A.M.Soc. 75 (1969) 1293-1298.

\bibitem{Qu73} D.G. Quillen. Higher Algebraic K-theory I. Lecture Notes in Math. 341 (1973) 85-147.

\bibitem{Sn79}  V.P. Snaith. {\em Algebraic Cobordism and K-theory}. Mem. Amer. Math. Soc. \#221 (1979).

\bibitem{Sn81}  V.P. Snaith.  Localised stable homotopy of some classifying spaces. Math. Proc. Camb. Phil. Soc. 89 (1981) 325-330.

\bibitem{Z} V.P. Snaith. A model for equivariant $K$-theory. C.R. Acad. Sci. Canada 3(3) (1981) 143--147.

\bibitem{Sn83a}  V.P. Snaith. Localised stable homotopy and algebraic K-theory. Memoirs Amer. Math.Soc. \#280 (1983).

\bibitem{Sn07} V.P. Snaith. {\em  Stable Homotopy Theory around the Arf-Kervaire invariant}. Preprint.

\bibitem{X} M. Spitzweck and P.A. {\O}stv{\ae}r.  The Bott inverted infinite projective space is homotopy algebraic $K$-theory.  Bulletin of the London Mathematical Society 41(2) (2009) 281-292.

\bibitem{TT} R. Thomason and T. Trobaugh.  Higher algebraic K-theory of schemes and derived categories. {\em The Grothendieck festschrift}, vol. 3 (1990), 247-436, Boston, Birkhauser.

\bibitem{Voev98} V. Voevodsky. ${\mathbb A}^{1}$-homotopy theory. Doc. Math. Extra Vol. ICM I (1998) 579-604.

\bibitem{Voev02b} V. Voevodsky. Open problems in the motivic stable homotopy theory I. {\em Motives, Polylogarithms and Hodge Theory} Part I (Irvine CA 1998) 3-34, Int. Press Lect. Ser. 3 I Int. Press Somerville MA (2002).

\bibitem{Voev04}  V. Voevodsky.  On the zero slice of the sphere spectrum. Proc. Steklov Inst. Math. (translation) (2004) no.3 (246) 93-102.

\bibitem{VSF00}  V. Voevodsky, A. Suslin and E.M. Friedlander. {\em Cycles, Transfers and Motivic Homology Theories}. Annals of Math. Studies \#143,
Princeton Univ. Press (2000).

\bibitem{Weibel88} C.A. Weibel. Homotopy algebraic K-theory. Contemp. Math. \#83 (1988) 461-488.


\end{thebibliography}
\end{document}